\documentclass{article}%
\usepackage{amssymb}
\usepackage{amsfonts}
\usepackage{graphicx}
\usepackage{amsmath}%
\providecommand{\U}[1]{\protect\rule{.1in}{.1in}}

\newenvironment{proof}[1][Proof]{\textbf{#1.} }{\ \rule{0.5em}{0.5em}}
\begin{document}
{\LARGE Invariant subspaces of generalized Hardy algebras associated with
compact abelian group actions on W*-algebras}

\bigskip

\ \ \ \ \ \ \ \ \ \ \ \ \ \ \ \ \ \ \ \ \ \ \ \ \ \ \ \ \ \ \ \ \ \ \ \ \ \ \ \ Costel
Peligrad

\bigskip

Department of Mathematical Sciences, University of Cincinnati, PO Box 210025,
Cincinnati, OH 45221-0025, USA. E-mail address: costel.peligrad@uc.edu

\bigskip


\bigskip


\bigskip

ABSTRACT. We consider an action of a compact abelian group whose dual is any
subgroup of the additive group of real numbers (so, an archimedean linearly
ordered group) or a direct product (or sum) of such groups on a W*-algebra,
$M$. We define the generalized Hardy subspace of the Hilbert space of a
standard representation the algebra, and the Hardy subalgebra of analytic
elements of $M$ with respect to the action. We find conditions in order that
the Hardy algebra is a hereditarily reflexive algebra of operators. In
particular if every non zero spectral subspace, contains a unitary operator,
the condition is satisfied and therefore the Hardy algebra is hereditarily
reflexive. This is the case if the action is the dual action on a crossed
product, or an ergodic action, or, if, in some situations, the fixed point
algebra is a factor.

\bigskip

\section{\bigskip Introduction}

This paper is concerned with the study of invariant subspaces and reflexivity
of operator algebras associated with compact group actions on W*-algebras.
Recall first the definition of a reflexive operator algebra.

Let $A\subset B(X)$ be a weakly closed algebra of operators on a Banach space
$X$. Denote by $Lat(A)$ the lattice of closed subspaces of $X$\ that are
invariant for all operators $a\in A.$ Let
\[
algLat(A)=\left\{  b\in B(X):bK\subset K\text{ for all }K\in Lat(A)\right\}
.
\]
The algebra $A$\ is called reflexive if $A=algLat(A)$. Hence, a reflexive
operator algebra is completely determined by the lattice of its invariant
subspaces. An algebra $A\subset B(X)$\ is called hereditarily reflexive if
every unital weakly closed subalgebra of $A$\ is reflexive. Sarason [19],
proved two results: (1) every commutative von Neumann algebra is hereditarily
reflexive and (2) the algebra of analytic Toeplitz operators on the Hardy
space $H^{2}(\boldsymbol{T})$ where $\boldsymbol{T}$\ is the the unit circle
$\boldsymbol{T}=\left\{  z\in%
\mathbb{C}
:\left\vert z\right\vert =1\right\}  ,$ is hereditarily reflexive. In [14] we
extended this result in two directions: (1) to the case of $H^{p}%
(\boldsymbol{T}),1<p<\infty$ and (2) to the not necessarily commutative case
of non selfadjoint crossed products of finite von Neumann algebras by the
semigroup $%
\mathbb{Z}
_{+}$.\ Later, in [9], Kakariadis has considered the more general case of
reduced w*-semicrossed products and, among other results, he has extended the
particular case of our reflexivity result in [14, Proposition 4.5,] for
$p=2$\ to the semicrossed product setting [9, 2.10.H.]. This result was
considered later by Helmer [5] in the context of W*-correspondences [12].
Further, in [15],\ we studied a related problem in a more general setting than
the crossed product or the reduced $w^{\ast}$semicrossed product considered in
[14] and [9, 2.9] for the case of von Neumann algebras. We considered a
W*-dynamical system $(M,\mathbf{T}\boldsymbol{,}\alpha)$ where $\mathbf{T=}%
\left\{  z\in%
\mathbb{C}
:\left\vert z\right\vert =1\right\}  \ $is the circle group$\ $and $M$\ is a
$\sigma-$finite W*-algebra. We constructed a standard covariant representation
of the system on a certain Hilbert space, $H$, a generalized Hardy space,
$H_{+}$\ and the corresponding Hardy algebra $M_{+}\subset B(H_{+}).$\ We have
shown that if the spectral subspace corresponding to the smallest positive
element of the spectrum contains a unitary element, then, the algebra $M_{+}$
is reflexive. Actually, [15, Theorem 3.5.] shows that if $M\subset B(H)\ $is a
$\sigma$-finite$\ $von Neumann algebra in its standard representation such
that each spectral subspace contains a unitary operator (as is, in particular,
the algebra of analytic Toeplitz operators considered by Sarason), then,
$M_{+}\subset B(H_{+})$\ is a reflexive operator algebra. Recently Bickerton
and Kakariadis [2] have obtained results about reflexivity of algebras
associated with actions of $%
\mathbb{Z}
_{+}^{d}$ (the direct product of $d$\ copies of $%
\mathbb{Z}
_{+}).$

In this paper we make two significant steps towards solving the reflexivity
problem of Hardy algebras associated to one-parameter dynamical systems
$(M,\mathbf{%
\mathbb{R}
}\boldsymbol{,}\alpha)$: 1. We consider a W*-dynamical system $(M,G,\alpha
)$\ where $M$\ is a von Neumann algebra in standard form, and 2. $G$\ is a
compact abelian group whose dual is an arbitrary subgroup of $%
\mathbb{R}
$\ (possibly $%
\mathbb{R}
$\ itself with the discrete topology), so, a discrete group with a linear
archimedean order. We also consider actions of compact abelian groups, $G,$
whose duals, $\Gamma,$ are direct products or direct sums of discrete groups
with linearly archimedean order and consider the lattice order on $\Gamma$
(see for instance [3]). In Section 2.1. we define a standard covariant
representation of the system that will be the framework for the rest of the
paper. In the Corollary to Proposition 2.4. we show that every von Neumann
algebra in standard form (in particular every maximal abelian von Neumann
algebra) is hereditarily reflexive, thus extending the first result of Sarason
mentioned above to every von Neumann algebra in its standard representation.
In Section 3 we consider the case when the dual $\Gamma$\ of $G$\ has an
archimedean linear order, or is a direct product of such groups, we define a
generalized Hardy space $H_{+}\subset H$ and a corresponding Hardy algebra
$M_{+}\subset B(H_{+}),$ where $H\ $is$.$the\ Hilbert space of the standard
covariant representation of the system $(M,G,\alpha)$ and we prove that, in
some conditions, including the conditions in [15], $M_{+}\subset B(H_{+})\ $is
hereditarily reflexive (in [15] we proved only reflexivity for the particular
case when $\Gamma=%
\mathbb{Z}
$)$.$We do not assume as in [15] that $M$\ is $\sigma$-finite. Also, if
$\Gamma$\ is an arbitrary archimedean linearly ordered discrete group, it can
be any subgroup of $%
\mathbb{R}
$\ with the discrete topology, not only $%
\mathbb{Z}
$ as in [15], [9], [14]. Examples include the Hardy algebra of analytic
Toeplitz operators, $H^{\infty}(\mathbf{T}),$ the results in [15], $w^{\ast}%
$-crossed products by abelian archimedean ordered discrete groups or a direct
product of such groups, some reduced $w^{\ast}$-semicrossed products
considered in [9], [2] and other situations as stated in the Corollaries
3.14., 3.15., 3.16. and 3.17.

\bigskip

\section{ Preliminary results and notations}

\subsection{\bigskip Standard representations of W*-algebras}

In this section we review some concepts and results related to the standard
representation of a von Neumann algebra. Some of these results are certainly
known, but we did not find an exact reference for them. We provide proofs of
these results for the convenience of the reader. In Proposition 2.4. and its
Corollary we prove that every von Neumann algebra in its standard
representation is hereditarily reflexive. In particular, every abelian von
Neumann algebra is hereditarily reflexive ([19, Theorem 1]).

\bigskip

Let $M$\ be a W*-algebra and let $\rho$\ be a weight on the positive part,
$M^{+},\ $of $M,$ that is a mapping $\rho:M^{+}\rightarrow\lbrack0,\infty
)\cup\left\{  \infty\right\}  $\ such that
\[
\rho(m+n)=\rho(m)+\rho(n),m,n\in M^{+}\
\]
and%
\[
\rho(\lambda m)=\lambda\rho(m),m\in M^{+},\lambda\in%
\mathbb{R}
,\lambda\geqslant0
\]
with the convention $0\cdot\infty=0.$\ As it is customary ([8], [20]), denote%
\[
\mathcal{N}_{\rho}=\left\{  m\in M:\rho(m^{\ast}m)<\infty\right\}  .
\]%
\[
N_{\rho}=\left\{  m\in M:\rho(m^{\ast}m)=0\right\}  .
\]%
\[
F_{\rho}=\left\{  m\in M^{+}:\rho(m)<\infty\right\}  .
\]%
\[
\mathcal{M}_{\rho}=linear\text{ }span\text{ }of\text{ }F_{\rho}.
\]
It is immediate that $\mathcal{N}_{\rho}$\ is a left ideal of $M$. The weight
$\rho$\ is called faithful if $N_{\rho}=\left\{  0\right\}  ,$\ normal if it
is the sum of a family $\left\{  \varphi_{\iota}\right\}  $ of positive normal
linear functionals and semifinite if $\mathcal{M}_{\rho}$\ or, equivalently
[20, 2.1.], $\mathcal{N}_{\rho}$\ \ is w*- dense in $M.$

Now let $M$ be a W*-algebra, $M_{0}\subset M$ a W*-subalgebra, and
$P_{0}:M\rightarrow M_{0}$\ a w*-continuous projection of norm $1$ of
$M$\ onto $M_{0}$ which is, in addition, faithful on the set of positive
elements of $M.$ Let $\rho_{0}$\ be a faithful normal semifinite weight on
$M$. It is known that such a weight exists. Indeed, consider a family
$\left\{  \varphi_{\iota}\right\}  $ of positive normal linear functionals of
$M_{0}\ $such that their supports $\left\{  p_{\iota}\right\}  $\ form a
maximal family of mutually orthogonal projections of $M_{0},$\ in particular
$\sum p_{\iota}=I.$ Then $\rho_{0}=\sum\varphi_{\iota}$ is a faithful normal
semifinite weight of $M_{0}.$

\bigskip

The following fact is stated in [19, Corollary 10.5] as a consequence of a
theorem of Takesaki [20, Theorem 10.1.]. We present a short proof of this fact
in our setting for the convenience of the reader.

\bigskip

\textbf{Lemma 2.1. }$\rho=\rho_{0}\circ P_{0}$ \textit{is a faithful normal
semifinite weight on} $M.$

\bigskip

\begin{proof}
Since $\rho_{0}$\ and $P_{0}$\ are faithful, it follows that $\rho$\ is
faithful. Since $\rho_{0}$\ is normal and $P_{0}$\ is w*-continuous, it
follows that $\rho$\ is normal. To prove that $\rho$\ is semifinite, notice
that from the definition of $\mathcal{N}_{\rho}$\ we have that $M\mathcal{N}%
_{\rho_{0}}\subset\mathcal{N}_{\rho},$ where $M\mathcal{N}_{\rho_{0}}%
$\ denotes the linear span of $\left\{  mn:m\in M,n\in\mathcal{N}_{\rho_{0}%
}\right\}  .$ Since $\mathcal{N}_{\rho_{0}}$ is w*-dense in $M_{0}$\ it
follows that $M\mathcal{N}_{\rho_{0}}$\ and therefore $\mathcal{N}_{\rho}$ is
dense in $M$\ so $\rho$\ is semifinite.
\end{proof}

\bigskip

By [8, Theorem 7.5.3.], there exists a faithful normal representation
$\pi_{\rho}$ of $M$ on the completion $H_{\rho}$ of $\mathcal{N}_{\rho}\subset
M$ with respect to the inner product
\[
\left\langle m,n\right\rangle =\rho(n^{\ast}m).
\]
This representation is uniquely determined up to unitary equivalence and is,
in that sense, independent of the choice of the weight $\rho_{0}.$ We will use
the version of Tomita-Takesaki Theorem from [8, Theorem 9.2.37.]. If $S$\ is
the conjugate linear operator defined on $\mathcal{N}_{\rho}\cap
\mathcal{N}_{\rho}^{\ast}$ by $S(n)=n^{\ast},$\ then $S$\ is a preclosed
densely defined operator on $H_{\rho}$\ and its closure has the polar
decomposition $J\Delta^{\frac{1}{2}},$ in which $\Delta$\ is an invertible
positive operator and $J$\ is a conjugate linear isometry acting on $H_{\rho}$
such that $J^{2}=I$ and $J\pi_{\rho}(M)J=\pi_{\rho}(M)^{^{\prime}}$\ where
$\pi_{\rho}(M)^{^{\prime}}$\ is the commutant of $\pi_{\rho}(M)$\ in
$B(H_{\rho}).$ We have $\mathcal{N}_{\rho_{0}}\subset\mathcal{N}_{\rho}$,
where, as above
\[
\mathcal{N}_{\rho_{0}}=\left\{  m\in M_{0}:\rho_{0}(m^{\ast}m)<\infty\right\}
.
\]
and%
\[
\mathcal{N}_{\rho}=\left\{  m\in M:\rho(m^{\ast}m)<\infty\right\}  .
\]
\bigskip

In the rest of the paper if $\rho$\ is a faithful, normal semifinite weight on
$M^{+}$\ we will identify $\pi_{\rho}(M)$ with $M$ and will write $m$\ instead
of $\pi_{\rho}(m),m\in M.$ We will call this representation the standard
representation of $M$\ and we will refer to the inclusion $M\subset B(H_{\rho
})$\ as the standard form of $M.$ Also, we will denote $H_{\rho}$\ by $H$ and
the closure of $\mathcal{N}_{\rho_{0}}$\ in $H$\ by $H_{0}.$

\bigskip

\textbf{Lemma 2.2. }\textit{The restriction of }$P_{0}$\textit{\ to
}$\mathcal{N}_{\rho}$\textit{\ extends to the orthogonal projection of }%
$H$\textit{\ onto }$H_{0}.$

\bigskip

\begin{proof}
Clearly, $P_{0}(\mathcal{N}_{\rho_{0}})=\mathcal{N}_{\rho_{0}}.$ We will prove
next that $P_{0}(\mathcal{N}_{\rho})\subset\mathcal{N}_{\rho_{0}}$. Indeed,
let $n\in\mathcal{N}_{\rho}.$ Then $n=P_{0}(n)+(I-P_{0})(n)=n_{0}+n_{1}.$
Since $n\in\mathcal{N}_{\rho},$ we have $\rho(n^{\ast}n)<\infty,$ so
\[
\rho(n_{0}^{\ast}n_{0}+n_{0}^{\ast}n_{1}+n_{1}^{\ast}n_{0}+n_{1}^{\ast}%
n_{1})=\rho_{0}(P_{0}(n_{0}^{\ast}n_{0}+n_{0}^{\ast}n_{1}+n_{1}^{\ast}%
n_{0}+n_{1}^{\ast}n_{1}))=
\]%
\[
\rho_{0}(n_{0}^{\ast}n_{0}+n_{0}^{\ast}P_{0}(n_{1})+P_{0}(n_{1}^{\ast}%
)n_{0}+P_{0}(n_{1}^{\ast}n_{1}))=
\]%
\[
\rho_{0}(n_{0}^{\ast}n_{0}+P_{0}(n_{1}^{\ast}n_{1}))<\infty.
\]
Therefore, $\rho_{0}(n_{0}^{\ast}n_{0})<\infty$ and we are done. On the other
hand,
\[
\left\langle P_{0}(n),m\right\rangle =\left\langle n,P_{0}(m)\right\rangle .
\]
since both of the above terms equal $\rho_{0}(P_{0}(m^{\ast})P_{0}(n))$, so
$P_{0}\ $is.self adjoint.
\end{proof}

\bigskip

\textbf{Lemma 2.3. }\textit{i) }$MH_{0}$\textit{\ is dense in }$H$%
\textit{\ and }$H_{0}$\textit{\ is a separating set for }$M$\textit{\ that is,
if }$m\in M$\textit{\ is such that }$m\xi_{0}=0$\textit{\ for all }$\xi_{0}\in
H_{0},$ \textit{then }$m=0.$ \textit{Here}, $MH_{0}$\ \textit{denotes the
linear span of} $\left\{  m\xi:m\in M,\xi\in H_{0}\right\}  .$

\textit{ii) }$M^{\prime}H_{0}$\textit{\ is dense in }$H$\textit{\ and }$H_{0}%
$\textit{\ is a separating set for }$M^{\prime}$\textit{\ that is, if
}$m^{\prime}\in M^{\prime}$\textit{\ is such that }$m^{\prime}\xi_{0}%
=0$\textit{\ for all }$\xi_{0}\in H_{0},$ \textit{then }$m^{\prime}=0.$

\bigskip

\begin{proof}
i) We will prove that $MH_{0}$\ is dense in $\mathcal{N}_{\rho}\subset
H$\ and, since $\mathcal{N}_{\rho}$\ is dense in $H,$\ the first part of i)
will follow. Let $x\in\mathcal{N}_{\rho}.$ Therefore, $\rho(x^{\ast}%
x)=\rho_{0}(P_{0}(x^{\ast}x))<\infty.$ Since $\rho_{0}$\ is a faithful normal
semifinite weight on $M_{0}$\ we can assume that $\rho_{0}$ is the sum of a
family of normal positive linear functionals $\left\{  \varphi_{\iota
}\right\}  $\ on $M_{0}$\ such that their suports $\left\{  p_{\iota}\right\}
$ form a maximal family of mutually ortogonal projections in $M_{0}$\ and
$\sum p_{\iota}=I.$ Since $\rho_{0}(P_{0}(x^{\ast}x))=\sum\varphi_{\iota
}(P_{0}(x^{\ast}x))<\infty,$ it follows that the summable family of positive
numbers $\left\{  \varphi_{\iota}(P_{0}(x^{\ast}x))\right\}  _{\iota}$ is at
most countable, say $\left\{  \varphi_{i}(P_{0}(x^{\ast}x))\right\}  _{i\in%
\mathbb{N}
}$ with $\rho(x^{\ast}x)=\rho_{0}(P_{0}(x^{\ast}x))=\sum_{i=1}^{\infty}%
\varphi_{i}(P_{0}(x^{\ast}x))<\infty.$ Let $q_{n}=\sum_{i=1}^{i=n}p_{i}%
,$\ where, for each $i\in%
\mathbb{N}
,$\ $p_{i}$\ is the supportof $\varphi_{i}.$ Then, $q_{n}\in\mathcal{N}%
_{\rho_{0}}\subset H_{0}$\ and $xq_{n}\in M\mathcal{N}_{\rho_{0}}\subset
MH_{0}$ for every $n\in%
\mathbb{N}
$. Clearly, since $p_{i}$\ is the support of $\varphi_{i}$\ we have
$\varphi_{i}(y)=\varphi_{i}(p_{i}y)=\varphi_{i}(yp_{i})$ for all $i\in%
\mathbb{N}
,y\in M_{0}$\ and\ $\varphi_{i}(q_{n}y)=\varphi_{i}(y)$\ for every $i,n\in%
\mathbb{N}
$\ with $i\leqslant n,y\in M_{0}$ and $\varphi_{i}(q_{n}x)=0$ if $i>n$\ We
will show that $\lim_{n\rightarrow\infty}xq_{n}=x$ in $H.$Let $\epsilon>0.$
Then, there exists $N=N_{\epsilon}>0$\ such that $\sum_{i=N+1}^{\infty}%
\varphi_{i}(P_{0}(x^{\ast}x))<\epsilon^{2}.$ Therefore, if $n\geqslant N$%
\[
\left\langle q_{n}x-x,q_{n}x-x\right\rangle =\rho((q_{n}x-x)^{\ast}%
(q_{n}x-x))=\sum_{i=1}^{n}\varphi_{i}(P_{0}((q_{n}x-x)^{\ast}(q_{n}x-x)))+
\]%
\[
+\sum_{i=n+1}^{\infty}\varphi_{i}(P_{0}((q_{n}x-x)^{\ast}(q_{n}x-x))
\]
Using the preceding observations, if $i\leqslant n,$ we get\
\[
\varphi_{i}(P_{0}((q_{n}x-x)^{\ast}(q_{n}x-x)))=
\]

\[
\varphi_{i}(q_{n}P_{0}(x^{\ast}x)q_{n}-q_{n}P_{0}(x^{\ast}x)-P_{0}(x^{\ast
}x)q_{n}+P_{0}(x^{\ast}x))=0
\]
and, if $i>n\geqslant N$
\[
\varphi_{i}(P_{0}((q_{n}x-x)^{\ast}(q_{n}x-x)=\varphi_{i}(P_{0}(x^{\ast}x))
\]
Hence $\left\langle q_{n}x-x,q_{n}x-x\right\rangle <\epsilon^{2}$ and the
claim is proven. To prove the second part, let $m\in M$ such that $m\xi_{0}%
=0$\ for every $\xi_{0}\in H_{0},$\ in particular, $mn=0$ for every
$n\in\mathcal{N}_{\rho_{0}}\subset H_{0}.$ Since $\mathcal{N}_{\rho_{0}}$ is
w*-dense in $M_{0}\subset M$ and $I\in M_{0},$ it follows that $m=0$\ and part
i) is proven.

ii). Denote by $K$ the closure of $M^{\prime}H_{0}.$\ Then $K$\ is a closed
subspace of $H$\ which is invariant for every $m^{\prime}\in M^{\prime},$
so,\ the orthogonal projection, $p,$ of $H$\ on $K$\ comutes with $M^{\prime
},$\ and therefore $p\in M.$ Since $H_{0}\subset K,$\ it follows that
$(1-p)H_{0}=\left\{  0\right\}  .$ Since by i) $H_{0}$\ is a separating set
for $M,$\ we have $1-p=0$ and thus $K=H.$To prove the second part of ii), let
$m^{\prime}\in M^{\prime}$ be such that $m^{\prime}H_{0}=\left\{  0\right\}
.$ It follows that $Mm^{\prime}H_{0}=\left\{  0\right\}  ,$ so $m^{\prime
}MH_{0}=\left\{  0\right\}  .$\ Since, by i) $MH_{0}$\ is dense in $H$, it
follows that $m^{\prime}=0.$
\end{proof}

\bigskip

Some of the statements in the next Proposition are probably known, but we did
not find a reference for any of them.

\bigskip

\textbf{Proposition 2.4.} \textit{i) Let }$M\subset B(H)$\textit{\ be a von
Neumann algebra in standard form. Then, every normal linear functional,
}$\varphi,$\textit{ on }$M$\textit{\ is a vector functional, that is, there
exist }$\xi,\eta\in H$\textit{ such that }$\varphi(m)=\left\langle m\xi
,\eta\right\rangle ,m\in M.$

\textit{ii) If }$N\subset B(H)$ \textit{is an abelian von Neumann algebra, not
necessarily in standard form, then every normal linear functional on }%
$N$\textit{\ is a vector functional.}

\textit{iii) If }$M_{0}\subset B(H_{0})$ \textit{is a maximal abelian}
\textit{von Neumann algebra, then it is spatially isomorphic with its standard
form.}

\bigskip

\begin{proof}
i) Let $p\in M$\ be a countably decomposable projection. According to [8,
9.6.18.], the hypotheses of [8, 9.6.20.] are satisfied, so\ there exists
$\xi_{0}\in H$ such that $J\xi_{0}=\xi_{0}$ and $\overline{M^{\prime}\xi_{0}%
}=pH.$ Therefore, every countably decomposable projection $p\in M\ $is a
cyclic projection. Now, let $\varphi$\ be a normal linear functional on $M.$
By the polar decomposition of normal linear functionals [8, Theorem 7.3.2.],
it is enough to prove the statement for normal positive linear functionals.
Let $\psi$\ be a normal positive functional on $M$\ and $p$\ its support (that
is, $p$ is the complement of the supremum of all projections $q\in M$\ for
which $\psi(q)=0).$ Then, $p$\ is countably decomposable, so by the previous
arguments, $p$ \ is a cyclic projection. Applying [8, Proposition 7.2.7.] it
follows that $\psi$\ is a vector normal positive functional.

ii) Let $\varphi$\ be a normal linear functional on $M.$\ As argued in i),
using the polar decomposition of normal linear functionals it is enough to
prove the statement in ii) for normal positive functionals. Let $\psi$\ be a
normal positive functional on $M$\ and $p$\ its support which is a countably
decomposable projection.\ Without loss of generality we can assume that $p=I.$
We will show that there exists a separating vector, $\xi_{0}\in H$ for $N$ and
therefore, cyclic for $N^{\prime}.$\ Let $\left\{  \xi_{i}\in H:i\in
A\right\}  $ be a maximal family of orthogonal unit vectors such that, the
projections $\left\{  p_{i}\in N:i\in A\right\}  $ onto $\left\{
H_{i}=\overline{N^{\prime}\xi_{i}}\subset H:i\in A\right\}  $ are mutually
orthogonal, so $H=\sum^{\oplus}H_{i}$. Since $p=I$\ is countably decomposable,
it follows that the set $A$\ is at most countable. Suppose $A\subseteq%
\mathbb{N}
.$ Let $\xi_{0}=\sum_{i\in A}\frac{1}{2^{i}}\xi_{i}.$ and let $m\in N^{+}$ be
such that $m\xi_{0}=0.$ Hence $\left\langle m\xi_{0},\xi_{i}\right\rangle
=0,i\in A.$\ Then, since $N$\ is abelian, so $N\subseteq N^{\prime},$ it
follows that $\left\langle m\xi_{i},\xi_{i}\right\rangle =0$ for every $i.$
Therefore, since $m\in N^{+},$ it follows that$\ m\xi_{i}=0,$ so
$mH_{i}=\left\{  0\right\}  $ for every $i,$ and thus $m=0.$ Hence $\xi_{0}$
is separating for $N.$ The statement ii) follows from [8, 7.2.7.].

iii) Let $\left\{  p_{\iota}\right\}  $ be a maximal family of mutually
orthogonal countably decomposable projections of $M_{0}$\ and $\left\{
\varphi_{\iota}\right\}  $ a family of positive linear functionals such that
the support of $\varphi_{\iota}$\ is $p_{\iota}.$ Clearly $\sum p_{\iota}=I.$
Let $\rho_{0}=\sum\varphi_{\iota}.$ Then $\rho_{0}$\ is a faithful normal
semifinite weight on $M_{0}^{+}.$ By ii) for every $\iota$\ there exists
$\xi_{\iota}\in p_{\iota}H_{0}$ such that $\varphi_{\iota}(m)=\left\langle
m\xi_{\iota},\xi_{\iota}\right\rangle ,m\in M_{0}.$ Obviously, $\xi_{\iota}$
is a cyclic and separating vector of $Mp_{\iota}|_{p_{\iota}H_{0}}$ for every
$\iota.$It is also clear that $\left\langle m_{1}\xi_{\iota},m_{2}\xi_{\iota
}\right\rangle =\varphi_{\iota}(m_{2}^{\ast}m_{1}),m_{1},m_{2}\in M_{0},$ for
every $\iota.$ If $\mathcal{N}_{\rho_{0}}$\ is as above,%
\[
\mathcal{N}_{\rho_{0}}=\left\{  m\in M_{0}:\rho_{0}(m^{\ast}m)<\infty\right\}
,
\]
then the mapping $mp_{\iota}\rightarrow m\xi_{\iota}$ extends to a unitary
operator from $H_{\rho_{0}}$ to $H_{0}$ and we are done.
\end{proof}

\bigskip

\textbf{Corollary }\textit{i) Every von Neumann algebra in standard form is
hereditarily reflexive.}

\textit{ii) [19, Theorem 2] Every abelian von Neumann algebra, not necessarily
in standard form is hereditarily reflexive.}

\bigskip

\begin{proof}
i) Follows from Proposition 2.4. i) and [10, Theorem 3.5.].

ii) Follows from Proposition 2.4. ii) and [10 Theorem 3.5.].
\end{proof}

\subsection{\bigskip W*-dynamical systems with compact abelian groups}

\bigskip

Let $(M,G\boldsymbol{,}\alpha)$ be a W*-dynamical system, where $M$ is a
$W^{\ast}-$algebra, $G$\ is a compact abelian group with dual $\Gamma$,\ and
$\alpha\ $a faithful $w^{\ast}-$continuous action of $G$\ on $M,.$that is
$\alpha_{g}\neq id$\ if $g\neq0,$ where $id$\ is the identity automorphism of
$M$\ and the mapping $g\rightarrow\varphi(\alpha_{g}(m))$\ for every $m\in M$
and every\ $\varphi\in M_{\ast},$\ where $M_{\ast}$\ denotes the predual of
$M.$ For each $\gamma\in\Gamma,$\ denote by
\[
M_{\gamma}=\left\{
{\displaystyle\int}
\overline{\left\langle g,\gamma\right\rangle }\alpha_{g}(m)dg:m\in M\right\}
.
\]
where the integral is taken in the $w^{\ast}-$topology.\ In particular, if
$\gamma=0,$\ $M_{0}$\ is the fixed point algebra of the system. It can
immediately be checked that
\[
M_{\gamma}=\left\{  m\in M:\alpha_{g}(m)=\left\langle g,\gamma\right\rangle
m\right\}
\]
It is clear that the mapping $P_{\gamma}:M\rightarrow M_{\gamma}$ defined by
$P_{\gamma}(m)=%
{\displaystyle\int}
\overline{\left\langle g,\gamma\right\rangle }\alpha_{g}(m)dg$ is a
w*-continuous projection of $M$\ onto the closed subspace $M_{\gamma}\subset
M.$ In particular, $P_{0}$\ is a w*-continuous projection of $M$\ onto $M_{0}%
$\ which is clearly faithful (on $M^{+}$). It is well known that $M$\ is the
$w^{\ast}$-closed linear span of $\left\{  M_{\gamma}:\gamma\in\Gamma\right\}
.$ The Arveson spectrum of the action $\alpha$\ is, by definition ([1], [13])
\[
sp(\alpha)=\left\{  \gamma\in\Gamma:M_{\gamma}\neq\left\{  0\right\}
\right\}  .
\]

\bigskip

\textbf{Lemma 2.5. }\textit{i) }$M_{-\gamma}=M_{\gamma}^{\ast},$where
$M_{\gamma}^{\ast}=\left\{  m^{\ast}:m\in M_{\gamma}\right\}  $ $M_{\gamma
}^{\ast}=\left\{  m^{\ast}:m\in M_{\gamma}\right\}  $.

\textit{ii) }$M_{\gamma_{1}}M_{\gamma_{2}}\subset M_{\gamma_{1}+\gamma_{2}}$
\textit{where} $M_{\gamma_{1}}M_{\gamma_{2}}$ \textit{is the linear span of
}$\left\{  xy:x\in M_{\gamma_{1}},y\in M_{\gamma_{2}}\right\}  .$

\textit{iii) If }$m\in M_{\gamma}$\textit{ has polar decomposition
}$m=u\left\vert x\right\vert ,$\textit{ then }$u\in M_{\gamma}$\textit{ and
}$\left\vert x\right\vert \in M_{0}.$

\bigskip

\begin{proof}
i) and ii) are obvious. iii) is a straightforward consequence of the
uniqueness of the polar decomposition of $m$.
\end{proof}

\bigskip

Let $(M,G,\alpha)$\ be as above, $\rho_{0}$\ a faithful normal semifinite
weight on $M_{0}$ and $\rho=\rho_{0}\circ P_{0}.$ Consider the corresponding
normal faithful representation $\pi_{\rho}$ on $H_{\rho}$ and the
Tomita-Takesaki operators $S,J$\ as in 2.1. above$.$ As in 2.1. we will write
$H$\ instead of $H_{\rho}$\ and $M$\ instead of $\pi_{\rho}(M).$ In the case
when $\Gamma$\ is a partially ordered group, this representation will allow us
to construct a generalized Hardy space on which, in certain situations, the
subalgebra of analytic elements of the system $(M,G\boldsymbol{,}\alpha)$\ is
hereditarily reflexive.

\bigskip

For every $g\in G$\ define the unitary operator $U_{g}\in B(H)$\ as the unique
extension of $U_{g}(n)=\alpha_{g}(n),n\in\mathcal{N}_{\rho}$ to $H.$ Then,
since clearly, $\mathcal{N}_{\rho}$ is an $\alpha$-invariant\ left ideal of
$M,$\ it is straightforward to check that the group of unitary operators
$\left\{  U_{g}:g\in G\right\}  $ implements the action $\alpha.$ Also, from
the definition of $S$ it follows that $SU_{g}=U_{g}S$\ and $S^{\ast}%
U_{g}=U_{g}S^{\ast}$\ for all $g\in G.$ Therefore, $JU_{g}=U_{g}J$, $g\in G.$
It follows that the group $\left\{  U_{g}:g\in G\right\}  $\ implements an
action $\alpha^{\prime}$\ of $G$\ on $M^{^{\prime}},$ namely%
\[
\alpha_{g}^{\prime}(JmJ)=U_{g}JmJU_{g}^{\ast}=J\alpha_{g}(m)J.
\]
Similarly with the projections $P_{\gamma}$\ of $M$\ onto $M_{\gamma}%
,\gamma\in\Gamma$\ one can define the projections $P_{\gamma}^{\prime}$\ of
$M^{\prime}$\ onto $M_{\gamma}^{\prime}=\left\{  x\in M^{\prime}:\alpha
_{g}^{\prime}(x)=\left\langle g,\gamma\right\rangle x\right\}  $%
\[
P_{\gamma}^{\prime}(x)=\int\overline{\left\langle g,\gamma\right\rangle
}\alpha_{g}^{\prime}(x)dg.
\]

\bigskip

The proof of the following lemma is a straghtforward application of the definitions.

\bigskip

\textbf{Lemma 2.6. }\textit{With the notations above, we have the following:}

\textit{i) }$\alpha_{g}^{\prime}(m^{\prime})=U_{g}m^{\prime}U_{g}^{\ast}%
$\textit{ is an action of }$G$\textit{\ on }$M^{\prime},$\textit{ where
}$M^{\prime}$\textit{\ is the commutant of }$M$\textit{\ in B(H).}

\textit{ii) If }$g\in G,$\textit{ then}$\ U_{g}$\textit{\ commutes with }%
$J$\textit{, and }$J\alpha_{g}(m)J=\alpha_{g}^{\prime}(JmJ),m\in M,g\in G.$

\textit{iii) (}$M^{\prime})_{\gamma}=JM_{-\gamma}J,\gamma\in\Gamma.$

\textit{iv) }$U_{g}(m^{\prime}\xi)=\alpha_{g}^{\prime}(m^{\prime})\xi,\xi\in
H_{0},m^{\prime}\in M^{\prime}.$

\textit{v) }$sp(\alpha)=sp(\alpha^{\prime}).$

\bigskip

We will need also the following

\bigskip

\textbf{Remark 2.7. }\textit{ If }$M$\textit{\ is a finite W*-algebra, then
}$M^{\prime}$\textit{\ is a finite W*-algebra.This fact is immediate from the
definition of the standard representations.}

\bigskip

Let $H_{\gamma}=\left\{
{\displaystyle\int}
\overline{\left\langle g,\gamma\right\rangle }U_{g}(\xi)dg:\xi\in H\right\}
=\left\{  \xi\in H:U_{g}\xi=\left\langle g,\gamma\right\rangle \xi\right\}  .$
Then, the map $P_{\gamma}^{H}$ from $H$ to $H_{\gamma}$ defined as follows%
\[
P_{\gamma}^{H}(\xi)=%
{\displaystyle\int}
\overline{\left\langle g,\gamma\right\rangle }U_{g}(\xi)dg:\gamma\in\Gamma
,\xi\in H.
\]
is an orthogonal projection of $H$\ onto the closed supspace $H_{\gamma}.$
Applying Lemma 2.2., we see that if $\gamma=0,$\ the Hilbert subspace
$H_{0}\subset H$\ coincides with the Hilbert subspace $H_{0}$\ considered in
Section 2.1.

\bigskip

\textbf{Lemma 2.8. }\textit{i) If }$\gamma_{1}\neq\gamma_{2},$\textit{ then
}$H_{\gamma_{1}}$\textit{ and }$H_{\gamma_{2}}$\textit{ are orthogonal.}

\textit{ii) For every }$\gamma\in sp(\alpha),$\textit{we have }$\overline
{M_{\gamma}H_{0}}=H_{\gamma},$\textit{ where }

$M_{\gamma}H_{0}=\left\{  m\xi_{0}:m\in M_{\gamma},\xi_{0}\in H_{0}\right\}
.$

\textit{iii) The direct sum of Hilbert spaces }$%
{\displaystyle\sum}
H_{\gamma}$\textit{ equals }$H.$

\textit{iv) For every }$\gamma\in sp(\alpha)$\textit{ we have }$\overline
{(M^{^{\prime}})_{\gamma}H_{0}}=H_{\gamma},$\textit{ where }

\textit{(}$M^{^{\prime}})_{\gamma}H_{0}=\{m^{^{\prime}}\xi_{0}:m^{^{\prime}%
}\in(M^{^{\prime}})_{\gamma},\xi_{0}\in H_{0}\}.$

\textit{v) For all }$\gamma,\gamma^{\prime}\in sp(\alpha)$ we have\textit{\ }%
$M_{\gamma}H_{\gamma^{\prime}}\subset H_{\gamma+\gamma^{\prime}}$ and
$M_{\gamma}^{\prime}H_{\gamma^{\prime}}\subset H_{\gamma+\gamma^{\prime}}.$

\bigskip

\begin{proof}
i) Let $\xi\in H_{\gamma_{1}},\eta\in H_{\gamma_{2}}.$ Then, by definition,
$U_{g}(\xi)=\left\langle g,\gamma_{1}\right\rangle \xi$ and $U_{g}%
(\eta)=\left\langle g,\gamma_{2}\right\rangle \eta,$ for all $g\in G.$ Since
the operators $U_{g}$ are unitary, we have%
\[
\left\langle \xi,\eta\right\rangle =\left\langle U_{g}\xi,U_{g}\eta
\right\rangle =\left\langle g,\gamma_{1}-\gamma_{2}\right\rangle \left\langle
\xi,\eta\right\rangle ,g\in G.
\]
Hence, if $\gamma_{1}\neq\gamma_{2}$ it follows that $\left\langle \xi
,\eta\right\rangle =0.$

ii) Since, by Lemma 2.3. i), $H_{0}$ is a cyclic set for $M,$ the subspace
$MH_{0}=\left\{  m\xi_{0}:m\in M,\xi_{0}\in H_{0}\right\}  $\ is dense in
$H.$\ Then, if $P_{\gamma}^{H}$ and $P_{\gamma}$ are the above projections, we
have
\[
M_{\gamma}H_{0}=P_{\gamma}(M)H_{0}=\left\{
{\displaystyle\int}
\overline{\left\langle g,\gamma\right\rangle }\alpha_{g}(m)\xi_{0}dg:m\in
M,\xi_{0}\in H_{0}\right\}  =
\]%
\begin{align*}
&  =\left\{
{\displaystyle\int}
\overline{\left\langle g,\gamma\right\rangle }U_{g}(\xi)dg:\xi=m\xi_{0},m\in
M,\xi_{0}\in H_{0}\right\}  =\\
&  .
\end{align*}%
\[
=\left\{  P_{\gamma}^{H}(\xi):\xi=m\xi_{0},m\in M,\xi_{0}\in H_{0}\right\}  .
\]
Since by Lemma 2.3. i) the subspace $MH_{0}$ is dense in $H$\ and $P_{\gamma
}^{H}$ is an orthogonal projection, the result stated in ii) follows.

iii) Let $\eta\in H$ be such that $\eta\perp H_{\gamma}$ for all $\gamma
\in\Gamma.$ Since $M$\ is the $w^{\ast}$-closed linear span of $\left\{
M_{\gamma}:\gamma\in\Gamma\right\}  ,$\ it follows from ii) that $\eta\perp
MH_{0},$ so, since by Lemma 2.3. i) $H_{0}$ is cyclic for $M,$ it follows that
$\eta=0.$

iv) The proof is similar with that of ii) taking into account that, according
to Lemma 2.3. ii), $H_{0}$ is cyclic for $M^{^{\prime}}$ as well.

v) Immediate from definitions.
\end{proof}

\bigskip

\section{Hereditary reflexivity of generalized Hardy algebras}

\bigskip

In this section we will construct the generalized Hardy space and the
generalized Hardy algebra and prove the main results of this paper, Theorem
3.8. and Theorem 3.9.

Let $\left(  M,G,\alpha\right)  $ be a W*-dynamical system with $G$ compact
abelian. Throughout this section we will assyme that $M\subset B(H)$ where $H$
is the Hilbert space constructed in Section 2.1. Suppose, in addition, that
$\Gamma$ is an archimedean linearly ordered discrete group, or $\Gamma
\subseteq\Pi_{\iota\in I}\Gamma_{\iota}$\ is the direct product or the direct
sum of archimedean linearly ordered (discrete) groups $\Gamma_{\iota}$. If
$\left(  \Gamma_{\iota}\right)  _{+}$ is the semigroup of non negative
elements of $\Gamma_{\iota},$ denote by $\Gamma_{+}=\Pi_{\iota\in I}\left(
\Gamma_{\iota}\right)  _{+}.$ Then, $\Gamma_{+}$ is a sub semigroup of
$\Gamma$\ such that
\[
\Gamma_{+}\cap(-\Gamma_{+})=\left\{  0\right\}
\]
and%
\[
\Gamma_{+}-\Gamma_{+}=\Gamma
\]
so $\Gamma_{+}$ defines a partial order on $\Gamma,$ namely $\gamma
_{1}\leqslant\gamma_{2}$ if $\gamma_{2}-\gamma_{1}\in\Gamma_{+}.$

\bigskip

\textbf{Lemma 3.1. }\textit{Let }$\gamma_{1},\gamma_{2}\in\Gamma_{+}%
,\gamma_{1}\neq\gamma_{2}.$\textit{ Then, either}

\textit{i) }$\gamma_{1},\gamma_{2}$\textit{ are not comparable under the above
order relation, or,}

\textit{ii) }$\gamma_{1}<\gamma_{2},$\textit{ or,}

\textit{iii) }$\gamma_{1}>\gamma_{2}$\textit{ and, in this case, there exists
}$p\in%
\mathbb{N}
$\textit{\ such that either }

\textit{iii a) }$\gamma_{1}\geqslant p\gamma_{2}$\textit{ and }$\gamma
_{1}<(p+1)\gamma_{2},$\textit{ or, }

\textit{iii b) }$\gamma_{1}\geqslant p\gamma_{2}$\textit{ and }$\gamma_{1}%
$\textit{ and }$(p+1)\gamma_{2}$\textit{ are not comparable.}

\bigskip

\begin{proof}
Suppose that $\gamma_{1},\gamma_{2}$\textit{ }are comparable. Thus, either
$\gamma_{1}<\gamma_{2},$ or $\gamma_{1}>\gamma_{2}.$ Suppose that $\gamma
_{1}>\gamma_{2}$. Since $\gamma_{1},\gamma_{2}\in\Pi_{\iota\in I}\left(
\Gamma_{\iota}\right)  _{+},$ we can write $\gamma_{1}=(\gamma_{1}^{\iota
})_{\iota\in I},\gamma_{2}=(\gamma_{2}^{\iota})_{\iota\in I}$ with $\gamma
_{1}^{\iota},\gamma_{2}^{\iota}\in\left(  \Gamma_{\iota}\right)  _{+},$ so
$\gamma_{1}^{\iota}\geqslant\gamma_{2}^{\iota},$ $\iota\in I$ and $\gamma
_{1}^{\iota_{0}}>\gamma_{2}^{\iota_{0}}$ for some $\iota_{0}\in I.$ Since for
every $\iota\in I$\ , $\Gamma_{\iota}$\ has an archimedean order, there exists
a largest $p_{\iota}\in%
\mathbb{N}
$\ such that $\gamma_{1}^{\iota}\geqslant p_{\iota}\gamma_{2}^{\iota}.$ If
$L\subset%
\mathbb{N}
$\ is the set of non repeating $p_{\iota}^{\prime}$s\ then, since $%
\mathbb{N}
$\ is well ordered, there exists $p=\min L.$ Therefore, $\gamma_{1}\geqslant
p\gamma_{2}.$ By the definition of $p\in%
\mathbb{N}
,$\ $\gamma_{1}\ngeqslant(p+1)\gamma_{2},$ so either iii a) or iii b) must hold.
\end{proof}

\bigskip

The following consequence of the above Lemma will be used

\bigskip

\textbf{Corollary 3.2. }\textit{Let }$\gamma_{0}\in\Gamma_{+}\smallsetminus
\left\{  0\right\}  $ and $\gamma\in\Gamma_{+}.$\textit{Then, there
exists\ }$p\in%
\mathbb{Z}
_{+}$\textit{\ such that either }$p\gamma_{0}\leqslant\gamma<(p+1)\gamma_{0}%
$\textit{\ or }$p\gamma_{0}\leqslant\gamma$\textit{\ and }$\gamma$\textit{\ is
not comparable with }$(p+1)\gamma_{0}.$

\bigskip

\begin{proof}
If $\gamma<\gamma_{0}$ or $\gamma$ and $\gamma_{0}$ are not comparable, then
$p=0$ satisfies the conclusion. If $\gamma>\gamma_{0},$\ the statement follows
from Lemma 3.1. iii).
\end{proof}

\bigskip

\textbf{Lemma 3.3. }\textit{If }$\Gamma\subseteq\Pi_{\iota\in I}\Gamma_{\iota
}$\textit{\ is the direct product or the direct sum of linearly ordered
discrete groups }$\Gamma_{\iota},$\textit{\ then }$\Gamma$\textit{ is lattice
ordered (see [3]), i.e. if }$A=\left\{  \gamma_{1},\gamma_{2},...\gamma
_{n}\right\}  $\textit{ is a finite subset of }$\Gamma$\textit{\ then there
exists }$\inf A$\textit{\ and }$\sup A$ in $\Gamma.$

\bigskip

\begin{proof}
If $\gamma_{j}=\left(  \gamma_{j}^{\iota}\right)  ,1\leqslant j\leqslant n,$
let $\mu^{\iota}=\min\left\{  \gamma_{j}^{\iota}:j=1,2,...n\right\}  $ and
$\nu^{\iota}=\max\left\{  \gamma_{j}^{\iota}:j=1,2,...n\right\}  $ for each
$\iota\in I.$ Then clearly $\inf A=\left(  \mu^{\iota}\right)  _{\iota\in I}$
and $\sup A=\left(  \nu^{\iota}\right)  _{\iota\in I}$.
\end{proof}

\bigskip

If $(M,G,\alpha),M\subset B(H)$, $H_{0},$ $\widehat{G}=\Gamma$ and $\Gamma
_{+}$\ are as above, define%
\[
H_{+}=\sum_{\gamma\in\Gamma_{+}}H_{\gamma}%
\]
Let $p_{+}$\ be the orthogonal projection of $H$\ onto $H_{+}$. By Lemma 2.8.
v), the (closed) subspace $H_{+}\subset H$\ is invariant for $\vee
_{\gamma\geqslant0}M_{\gamma},$ and for $\vee_{\gamma\geqslant0}(M^{\prime
})_{\gamma}.$We will denote by $M_{+}$\ the weak operator closure%
\[
M_{+}=\overline{p_{+}(\vee_{\gamma\in\Gamma_{+}}M_{\gamma})p_{+}}^{wo}%
\]
in $B(H_{+})$\ and similarly%
\[
(M^{\prime})_{+}=\overline{p_{+}(\vee_{\gamma\geqslant0}(M^{\prime})_{\gamma
})p_{+}}^{wo}%
\]
where $\vee_{\gamma\geqslant0}M_{\gamma},$ is the algebra generated by
$\left\{  M_{\gamma}:\gamma\geqslant0\right\}  )$ and $\vee_{\gamma\geqslant
0}(M^{\prime})_{\gamma}$ is the algebra generated by $\left\{  (M^{\prime
})_{\gamma}:\gamma\geqslant0\right\}  )$. Then, we will call $H_{+}$ the
generalized Hardy space and $M_{+}$\ the generalized Hardy algebra of analytic
elements of the dynamical system $(M,G,\alpha).$

\bigskip

To prove hereditary reflexivity, we also need the following

\bigskip

\textbf{Lemma 3.4. }\textit{If }$\psi$\textit{ is a weakly continuous
functional on }$M_{+}\subset B(H_{+}),$\textit{\ then }$\psi$\textit{\ is a
vector functional (that is, there exist }$\xi,\eta\in H_{+}$\textit{\ such
that }$\psi(m)=\left\langle m\xi,\eta\right\rangle $ for all $m\in M_{+}$).

\bigskip

\begin{proof}
Since $\psi$\ is weakly continuous, there exist $n\in%
\mathbb{N}
$ and $\xi_{i},\eta_{i}\in H_{+},1\leqslant i\leqslant n$ such that
$\psi(m_{+})=\sum_{i}\left\langle m_{+}\xi_{i},\eta_{i}\right\rangle ,m_{+}\in
M_{+}.$ Now let $\widetilde{\psi}$ be the functional on \ defined by
$\widetilde{\psi}(b)=\sum_{i}\left\langle b\xi_{i},\eta_{i}\right\rangle ,b\in
B(H).$ Since $\xi_{i},\eta_{i}\in H_{+},$ it follows that $\widetilde{\psi
}(b)=\widetilde{\psi}(p_{+}bp_{+}),b\in B(H),$ where, as above, $p_{+}$\ is
the projection of $H$\ onto $H_{+}.$ The restriction of $\widetilde{\psi}$\ to
$M$\ is a normal linear functional of $M.$ Applying Proposition 2.4. to this
restriction, it follows that there exist $\xi,\eta\in H$ such that
$\widetilde{\psi}(m)=\left\langle m\xi,\eta\right\rangle ,m\in M.$ Since, as
noticed before, $\widetilde{\psi}(m)=\widetilde{\psi}(p_{+}mp_{+}),$ we can
take $\xi,\eta\in H_{+}.$ Therefore, in particular, $\psi(m_{\gamma
})=\widetilde{\psi}(m_{\gamma})=\widetilde{\psi}(p_{+}mp_{+})=\psi(m_{\gamma
})=\left\langle m_{\gamma}\xi,\eta\right\rangle $ for every $m_{\gamma}\in
M_{\gamma},\gamma\in\Gamma_{+}.$ The definition of $M_{+}$\ implies that
$\psi(m_{+})=\left\langle m_{+}\xi,\eta\right\rangle ,m_{+}\in M_{+}$ and the
proof is completed.
\end{proof}

\bigskip

Loginov and \v{S}ul'man [10, Theorem 2.3.] have shown, in particular, that if
a reflexive algebra satisfies the hypothesis of Lemma 3.4. then it is
hereditarily reflexive, that is, all unital weakly closed subalgebras are
reflexive. Under the name of super reflexivity this fact has been also
considered in [4, Proposition 2.5. (1)].

\bigskip

In Theorem 3.8. below we will assume that $\left(  M,G,\alpha\right)  $
satisfies the following condition:

\bigskip

(\textbf{C}) For every $\gamma\in$ $sp(\alpha)\setminus\left\{  0\right\}  $
there exists an element $u_{\gamma}\in M_{\gamma}$ such that $u_{\gamma}%
^{\ast}u_{\gamma}=u_{\gamma}u_{\gamma}^{\ast}=e_{\gamma}$ where $e_{\gamma}$
is a central projection of $M$ and $M_{\gamma}$ = $M_{0}u_{\gamma}$. For
$\gamma=0$\ we will take $u_{0}=I.$

\bigskip

Examples of dynamical systems $\left(  M,G,\alpha\right)  $\ satisfying this
condition include the following:

\bigskip

a) If $M$\ is a finite W*-algebra and the center, $Z(M_{0}),$\ of $M_{0}$\ is
contained in the center, $Z(M),$\ of $M$ [18]$.$ This is the case, in
particular, when $M$\ is a finite W*-algebra and $M_{0}$\ is a factor (this
case will be discussed in a more general context in part b)). The conditions
$M$\ finite and $Z(M_{0})\subset Z(M)$ also hold if $M=\oplus M_{i}$\ and
$\left(  M_{i},G,\alpha\right)  $ is a finite W*-algebra and the fixed point
algebra is a factor. Corollary 3.14. below will refer to these example.

\bigskip

b) If $M$\ is a semifinite injective von Neumann algebra such that $M_{0}$\ is
a factor, except when $M$ is type $III$ and $M_{0}$ is a type $II_{1}$\ factor
[21]. Thomsen has proved that in these cases, the action $\alpha$\ has full
unitary spectrum, that is every nonzero spectral subspace contains unitary
operators. This is the case, in particular, when $M$\ is a finite W*-algebra
and $M_{0}$\ is a factor. In particular, this latter situation occurs if
$\alpha$\ is a prime action of the compact abelian group $G$\ on the
hyperfinite type $II_{1}$\ factor [6], [7], in particular if $\alpha$\ is
egodic. Recall that an action is called prime if the fixed point algebra is a
factor. In particular if the action \ $\alpha$ is faithful, then the all the
examples in this part b) satisfy $sp(\alpha)=\Gamma.$ Corollary 3.15. below
will refer to these examples. Also, the Condition (\textbf{C}) is satisfied if
$M$\ is the crossed product of a von Neumann algebra $M_{0.}$\ by an abelian
discrete group $\Gamma.$\ Corollaries 3.16. and 3.17. will consider this case.

\bigskip

\textbf{Lemma 3.5.\ }\textit{Suppose that condition (\textbf{C})}$.$\textit{is
satisfied. Then}

\textit{i) }$M_{\gamma}=u_{\gamma}M_{0}$

\textit{ii) There exists an element }$w_{\gamma}\in M^{\prime}=JMJ$%
\textit{\ such that }$w_{\gamma}w_{\gamma}^{\ast}=w_{\gamma}^{\ast}w_{\gamma
}=e_{\gamma},$\textit{ }and $(M^{\prime})_{\gamma}=(M^{\prime})_{0}w_{\gamma
}=w_{\gamma}(M^{\prime})_{0}$.

\bigskip

\begin{proof}
i) Clearly, if $x=mu_{\gamma}$\ for some \ then $x=me_{\gamma}u_{\gamma
}=e_{\gamma}mu_{\gamma}=u_{\gamma}(u_{\gamma}^{\ast}mu_{\gamma})\in u_{\gamma
}M_{0}$ and conversely.

ii) Obviously, $w_{\gamma}=Ju_{\gamma}^{\ast}J$ satisfies the equality
$w_{\gamma}w_{\gamma}^{\ast}=w_{\gamma}^{\ast}w_{\gamma}=Je_{\gamma}J$.Since
$e_{\gamma}$\ is a central projection of $M,$\ we can apply [8, 9.6.18.] to
get $Je_{\gamma}J=e_{\gamma}.$
\end{proof}

\bigskip

\textbf{Lemma 3.6.\ }\textit{Let }$(M,G,\alpha),M\subset B(H)$\textit{ be a
W*-dynamical system with }$G$\textit{\ compact abelian as above. Then, if
}$\gamma,\gamma^{\prime}\in sp(\alpha)$\textit{ and }$e_{\gamma}%
e_{\gamma^{\prime}}\neq0,$\textit{ we have }$\gamma^{\prime}-\gamma\in
sp(\alpha)$\textit{\ (therefore, }$\gamma-\gamma^{\prime}\in sp(\alpha
)$\textit{) and }$e_{\gamma}e_{\gamma^{\prime}}\leqslant e_{\gamma
-\gamma^{\prime}}.$

\bigskip

\begin{proof}
Since $e_{\gamma}e_{\gamma^{\prime}}\neq0,$ we have $u_{\gamma}u_{\gamma
}^{\ast}u_{\gamma^{\prime}}u_{\gamma^{\prime}}^{\ast}\neq0.$ so $u_{\gamma
}^{\ast}u_{\gamma^{\prime}}\neq0.$ Hence, applying Lemma 2.5. ii), it follows
that $M_{\gamma^{\prime}-\gamma}\neq0,$ so $\gamma^{\prime}-\gamma\in
sp(\alpha).$ To prove the last statement of the lemma, notice that by Lemma
2.5. and Lemma 3.5.
\[
e_{\gamma}e_{\gamma^{\prime}}=u_{\gamma}u_{\gamma}^{\ast}u_{\gamma^{\prime}%
}u_{\gamma^{\prime}}^{\ast}=u_{\gamma}u_{\gamma^{\prime}-\gamma}%
mu_{\gamma^{\prime}}^{\ast}=u_{\gamma}u_{\gamma^{\prime}-\gamma}%
me_{\gamma^{\prime}}u_{\gamma^{\prime}}^{\ast}=u_{\gamma}u_{\gamma^{\prime
}-\gamma}e_{\gamma_{^{\prime}}}mu_{\gamma^{\prime}}^{\ast}%
\]
for some $m\in M_{0}.$\ Further, using repeatedly Lemma 3.5. we get%
\[
e_{\gamma}e_{\gamma^{\prime}}=u_{\gamma}u_{\gamma^{\prime}-\gamma}%
e_{\gamma_{^{\prime}}}mu_{\gamma^{\prime}}^{\ast}=u_{\gamma}u_{\gamma^{\prime
}-\gamma}u_{\gamma^{\prime}}^{\ast}u_{\gamma^{\prime}}mu_{\gamma^{\prime}%
}^{\ast}=
\]%
\[
=u_{\gamma}u_{\gamma^{\prime}}^{\ast}(u_{\gamma^{\prime}}u_{\gamma^{\prime
}-\gamma}u_{\gamma^{\prime}}^{\ast})u_{\gamma^{\prime}}mu_{\gamma^{\prime}%
}^{\ast}=
\]%
\[
u_{\gamma-\gamma^{\prime}}m_{1}u_{\gamma^{\prime}-\gamma}m_{2}u_{\gamma
^{\prime}}mu_{\gamma^{\prime}}^{\ast}=u_{\gamma-\gamma^{\prime}}%
u_{\gamma^{\prime}-\gamma}m_{3}m_{2}e_{\gamma-\gamma^{\prime}}m_{1}%
u_{\gamma^{\prime}}mu_{\gamma^{\prime}}^{\ast}=e_{\gamma-\gamma^{\prime}}m_{4}%
\]
for some $m_{1},m_{2},m_{3},m_{4}\in M_{0}\backslash\left\{  0\right\}  .$
Therefore,
\[
e_{\gamma}e_{\gamma^{\prime}}=e_{\gamma-\gamma^{\prime}}m_{4}m_{4}^{\ast
}e_{\gamma-\gamma^{\prime}}\leqslant\left\Vert m_{4}m_{4}^{\ast}\right\Vert
e_{\gamma-\gamma^{\prime}}\
\]
So $e_{\gamma}e_{\gamma^{\prime}}\leqslant e_{\gamma-\gamma^{\prime}}.$
\end{proof}

\bigskip

\textbf{Lemma 3.7. }\textit{Suppose that Condition (\textbf{C}) is satisfied.
Then}

\textit{i) }$M_{+}$ \textit{is the w*-closed subalgebra of }$\mathit{B}%
$\textit{(}$H_{+}$\textit{) generated by }$M_{0}$\textit{\ and }

$\left\{  u_{\gamma}:\gamma\in sp(\alpha),\gamma>0\right\}  .$

\textit{ii) }$\left(  M^{\prime}\right)  _{+}$ \textit{is the w*-closed
subalgebra of }$\mathit{B}$\textit{(}$H_{+}$\textit{) generated by }$\left(
M^{\prime}\right)  _{0}$\textit{\ and }$\left\{  w_{\gamma}:\gamma\in
sp(\alpha^{\prime})=sp(\alpha),\gamma>0\right\}  .$

\bigskip

\begin{proof}
Follows from Lemma 3.5.
\end{proof}

\bigskip

We will prove next our results about reflexivity.

\bigskip

In Theorem 3.8. we assume that $\Gamma$\ is archimedean linearly ordered and
that Condition (\textbf{C}) is satisfied.

\bigskip

\textbf{Theorem 3.8. }\textit{Let }$(M,G,\alpha)$\textit{ be such that
}$\widehat{G}=\Gamma$ \textit{is archimedean linearly ordered and Condition
(\textbf{C}) is satisfied. Then }$M_{+}\subset B(H_{+})$\textit{ is
hereditarily reflexive}.

\bigskip

In Theorem 3.9. we assume that $\Gamma$ is a direct product (or a direct sum)
of archimedean linearly ordered discrete groups, but we assume a stronger
condition than Condition (\textbf{C}).

\bigskip\ 

\textbf{Theorem 3.9. }\textit{Let }$(M,G,\alpha)$\textit{ be such that
}$\widehat{G}=\Gamma\subseteq\Pi_{\iota\in I}\Gamma_{\iota}$ \textit{is the
direct product, or the direct sum}$,$ \textit{of archimedean linearly ordered
discrete groups, }$\Gamma_{\iota}.$ \textit{Suppose that }$sp(\alpha)=\Gamma
$\textit{\ and that for every }$\gamma\in sp(\alpha)=\Gamma,$\textit{\ there
exists a unitary operator }$u_{\gamma}\in M_{\gamma}.$ \textit{Then }%
$M_{+}\subset B(H_{+})$\textit{ is hereditarily reflexive.}

\bigskip

The proofs of these theorems will be given after some auxiliary results.

\bigskip

\textbf{Lemma 3.10. }\textit{Suppose that }$\Gamma$\textit{\ is archimedean
linearly ordered and that Condition (\textbf{C}) is satisfied. Then}

\textit{i) }$(M_{+})^{\prime}=(M^{\prime})_{+},$\textit{ where (}%
$M_{+})^{\prime}$\textit{ denotes the commutant of }$M_{+}$\textit{ in
}$B(H_{+})$ \textit{and}

\textit{ii) }$((M^{^{\prime}})_{+})^{^{\prime}}=M_{+}$

\bigskip

\begin{proof}
\textbf{ }i)\textbf{ }Let $x\in(M^{\prime})_{\gamma_{1}}$ and $m\in
M_{\gamma_{2}},\gamma_{1},\gamma_{2}\in sp(\alpha)\cap\Gamma_{+}$\ . Then,
\[
p_{+}xp_{+}mp_{+}=p_{+}xmp_{+}=p_{+}mxp_{+}=p_{+}mp_{+}xp_{+}%
\]
so, $(M^{\prime})_{+}\subset(M_{+})^{\prime}.$ To prove the converse
inclusion, let $x\in(M_{+})^{\prime}\subset B(H_{+}).$ Consider the following
dense subspace of $H_{+}$%
\[
H^{\prime}=\left\{  \sum_{\gamma\in F}H_{\gamma}:F\subset sp(\alpha)\cap
\Gamma_{+}\text{ a finite subset.}\right\}  .
\]
Then, clearly, the subspace
\[
H^{\prime\prime}=lin\left\{  u_{\gamma}\xi:\xi\in H^{\prime},\gamma\in
sp(\alpha)\right\}  .
\]
where $u_{\gamma}$\ is the partial isometry in Condition (\textbf{C}),\ is
dense in $H.$ Let $\eta=\sum_{i=1}^{n}u_{\gamma_{i}}\xi_{i}\in H".$\ Without
loss of generality, we will assume in the rest of this proof that $\gamma
_{1}\leqslant\gamma_{2}\leqslant...\leqslant\gamma_{n}.$ If $u_{\gamma_{i}%
}u_{\gamma_{i}}^{\ast}=u_{\gamma_{i}}^{\ast}u_{\gamma_{i}}=e_{\gamma_{i}%
},1\leqslant i\leqslant n$ are the central projections from Condition
(\textbf{C}) above, then, a standard calculation in the commutative W*-algebra
$Z(M),$ shows that%
\begin{equation}
e_{\gamma_{1}}\vee e_{\gamma_{2}}\vee...\vee e_{\gamma_{n}}=e_{\gamma_{1}%
}+\sum_{i=2}^{n}(1-e_{\gamma_{1}})...(1-e_{\gamma_{i-1}})e_{\gamma_{i}}.
\label{co1}%
\end{equation}
where $e_{\gamma_{1}}\vee e_{\gamma_{2}}\vee...\vee e_{\gamma_{n}}%
=\sup\left\{  e_{\gamma_{i}}:1\leqslant i\leqslant n\right\}  .$ Clearly, the
terms of the above sum are mutually orthogonal central projections in $Z(M).$
So if $\eta=\sum_{i=1}^{n}u_{\gamma_{i}}\xi_{i}\in H"$ it follows that%
\begin{equation}
\eta=e_{\gamma_{1}}\eta+\sum_{i=2}^{n}(1-e_{\gamma_{1}})...(1-e_{\gamma_{i-1}%
})e_{\gamma_{i}}\eta\text{.} \label{co2}%
\end{equation}
or,
\begin{equation}
\eta=\sum_{i=1}^{n}p_{i}\eta\label{co3}%
\end{equation}
where $p_{1}=e_{\gamma_{1}}$\ and $p_{i}=(1-e_{\gamma_{1}})...(1-e_{\gamma
_{i-1}})e_{\gamma_{i}},2\leqslant i\leqslant n.$ Define the operator
$\widehat{x}$\ on $H^{\prime\prime}$\ as follows%
\[
\widehat{x}(\sum_{i=1}^{n}u_{\gamma_{i}}\xi_{i})=\sum_{i=1}^{n}u_{\gamma_{i}%
}x\xi_{i}.
\]
We prove first that $\widehat{x}$\ is well defined. Indeed, suppose that
$\sum_{i=1}^{n}u_{\gamma_{i}}\xi_{i}=0.$ We will show that $\sum_{i=1}%
^{n}u_{\gamma_{i}}x\xi_{i}=0.$Now, if $\eta=\sum_{i=1}^{n}u_{\gamma_{i}}%
\xi_{i}=0,$ it follows that
\[
p_{i}\eta=0,1\leqslant i\leqslant n.
\]
We must show that
\begin{equation}
p_{i}\sum_{j=1}^{n}u_{\gamma_{j}}x\xi_{j}=0,1\leqslant i\leqslant n.
\label{co4}%
\end{equation}
These equalities imply that $\sum_{i=1}^{n}u_{\gamma_{i}}x\xi_{i}=0,$ so
$\widehat{x}$ is well defined. Since $p_{i}\eta=0,$ we have%

\[
\sum_{j=1}^{n}p_{i}u_{\gamma_{j}}\xi_{j}=\sum_{j=i}^{n}p_{i}u_{\gamma_{j}%
}p_{i}\xi_{j}=0.
\]
Thus, factoring out $u_{\gamma_{i}}$
\[
u_{\gamma_{i}}\sum_{j=1}^{n}p_{i}u_{\gamma_{i}}^{\ast}u_{\gamma_{j}}p_{i}%
\xi_{j}=0
\]
By multiplying the above equality by $u_{\gamma_{i}}^{\ast},$ and taking into
account that $p_{i}\leqslant e_{\gamma_{i}}$\ we get
\[
\sum_{j=i}^{n}p_{i}u_{\gamma_{i}}^{\ast}u_{\gamma_{j}}\xi_{j}=0.
\]
so%
\[
p_{i}e_{\gamma_{i}}\xi_{i}+\sum_{j=i+1}^{n}p_{i}u_{\gamma_{i}}^{\ast}%
u_{\gamma_{j}}\xi_{j}=0.
\]
where $u_{\gamma_{i}}^{\ast}u_{\gamma_{j}}\in M_{\gamma_{j}-\gamma_{i}}%
.$\ Since $\gamma_{i}\leqslant\gamma_{j},$ if $i\leqslant j\leqslant n$\ so,
$M_{\gamma_{j}-\gamma_{i}}\subset M_{+}$\ and $x\in(M_{+})^{\prime},$ it
follows that%
\[
x(p_{i}e_{\gamma_{i}}\xi_{i}+\sum_{i=2}^{n}p_{i}u_{\gamma_{i}}^{\ast}%
u_{\gamma_{j}}\xi_{j})=p_{i}e_{\gamma_{i}}x\xi_{i}+\sum_{j=i+1}^{n}%
p_{i}u_{\gamma i}^{\ast}u_{\gamma_{j}}x\xi_{j}=0.
\]
By multiplying the above equality by $u_{\gamma_{i}},$\ we get%
\[
p_{i}u_{\gamma_{i}}x\xi_{i}+\sum_{j=i+1}^{n}p_{i}u_{\gamma_{j}}x\xi_{j}%
=p_{i}\sum_{j=i}^{n}u_{\gamma_{j}}x\xi_{j}=p_{i}\sum_{j=i}^{n}u_{\gamma_{j}%
}x\xi_{j}=0.
\]
and this proves (\ref{co4}).\ Therefore, $\widehat{x}$\ is well defined. From
the definition of $\widehat{x}\ $it follows that $\widehat{x}m=m\widehat{x}%
$\ on $H"$ for all $m\in M$ and $\widehat{x}(\eta)=x(\eta)$ for every $\eta\in
H^{\prime}$, so if, as we will prove, $\widehat{x}\ $is bounded, it follows
that $\widehat{x}\in M^{\prime}.$\ Next we prove that the operator
$\widehat{x}$\ is bounded. Indeed, if, as above, $\eta=\sum u_{\gamma_{i}%
}e_{\gamma_{i}}\xi_{i},$ then, using the equality (\ref{co3}) and the fact
that $\widehat{x}$ commutes with $p_{i},1\leqslant i\leqslant n,$ we have
\[
\left\Vert \widehat{x}(\eta)\right\Vert =\left\Vert \widehat{x}\left(
\sum_{i=1}^{n}p_{i}\eta\right)  \right\Vert =\left\Vert \sum_{i=1}^{n}%
p_{i}\widehat{x}\left(  p_{i}\eta\right)  \right\Vert =\sqrt{\sum_{i=1}%
^{n}\left\Vert p_{i}\widehat{x}\left(  p_{i}\eta\right)  \right\Vert ^{2}}.
\]
Further, since $\gamma_{i}\leqslant\gamma_{j},$ and $p_{j}u_{\gamma_{i}}%
=0$\ when $i\leqslant j$\ and $x\in\left(  M_{+}\right)  ^{\prime}$ we have%
\[
\left\Vert p_{i}\widehat{x}\left(  p_{i}\eta\right)  \right\Vert
^{2}=\left\Vert p_{i}\widehat{x}(\sum_{j=1}^{n}p_{i}u_{\gamma_{j}}\xi
_{j})\right\Vert ^{2}=\left\Vert p_{i}\sum_{j=i}^{n}p_{i}u_{\gamma_{j}}%
x\xi_{j}\right\Vert ^{2}\leqslant
\]%
\[
\leqslant\left\Vert p_{i}\sum_{j=i}^{n}p_{i}u_{\gamma_{j}}x\xi_{j}\right\Vert
^{2}=\left\Vert p_{i}u_{\gamma_{i}}\sum_{j=i}^{n}p_{i}u_{\gamma_{i}}^{\ast
}u_{\gamma_{j}}x\xi_{j}\right\Vert ^{2}\leqslant\left\Vert p_{i}\sum_{j=i}%
^{n}p_{i}u_{\gamma_{i}}^{\ast}u_{\gamma_{j}}x\xi_{j}\right\Vert ^{2}=
\]%
\[
=\left\Vert p_{i}x\sum_{j=i}^{n}p_{i}u_{\gamma_{i}}^{\ast}u_{\gamma_{j}}%
x\xi_{j}\right\Vert ^{2}\leqslant\left\Vert x\right\Vert ^{2}\left\Vert
p_{i}\sum_{j=1}^{n}p_{i}u_{\gamma_{i}}^{\ast}u_{\gamma_{j}}\xi_{j}\right\Vert
^{2}=
\]%
\[
=\left\Vert x\right\Vert ^{2}\left\Vert p_{i}u_{\gamma_{i}}^{\ast}%
\eta\right\Vert ^{2}\leqslant\left\Vert x\right\Vert ^{2}\left\Vert p_{i}%
\eta\right\Vert ^{2}.
\]
Therefore, $\left\Vert \widehat{x}(\eta)\right\Vert \leqslant$ $\left\Vert
x\right\Vert \left\Vert \eta\right\Vert $ so $\widehat{x}$\ is bounded. As
noticed above, $\widehat{x}\in M^{\prime}$\ and, since obviously,
$p_{+}\widehat{x}p_{+}=x$ it follows that $x\in(M^{\prime})_{+}$ and we are done.

ii) follows from i) by replacing $M$\ with $M^{^{\prime}}.$
\end{proof}

\bigskip

The following version of Lemma 3.10. will be used in the proof of Theorem 3.9.
Here and in Theorem 3.9 we assume that $\Gamma$ is a direct product or a
direct sum of archimedean linearly ordered discrete groups, but we also assume
a stronger version of Condition (\textbf{C}), namely that for every $\gamma\in
sp(\alpha),M_{\gamma}$ contains a unitary operator and that the action is
faithful (i.e. $\alpha_{g}=id$\ implies $g=0$). These conditions imply that
$sp(\alpha)=\Gamma$ (see for instance [21 Lemma 2.2.])$.$

\bigskip

\textbf{Lemma 3.11.} \textit{Let }$(M,G,\alpha)$\textit{ be such that
}$\widehat{G}=\Gamma$ \textit{is a direct product, }$\Gamma=\Pi_{\iota\in
I}\Gamma_{\iota},$ (or a direct sum) \textit{of archimedean linearly ordered
discrete groups }$\Gamma_{\iota}.$ \textit{Suppose that }$\alpha$\textit{\ is
faithful and for every }$\gamma\in sp(\alpha),$\textit{\ there exists a
unitary operator }$u_{\gamma}\in M_{\gamma}$ (\textit{as noticed above, these
conditions imply that} $sp(\alpha)=\Gamma$). \textit{Then}

\textit{i) }$(M_{+})^{\prime}=(M^{\prime})_{+},$\textit{ where (}%
$M_{+})^{\prime}$\textit{ denotes the commutant of }$M_{+}$\textit{ in
}$B(H_{+})$ \textit{and}

\textit{ii) }$((M^{^{\prime}})_{+})^{^{\prime}}=M_{+}$.

\bigskip

\begin{proof}
i) As in the proof of the previous Lemma 3.10. it follows that $(M^{\prime
})_{+}\subset(M_{+})^{\prime}.$ To prove the opposite inclusion, let
$x\in(M_{+})^{\prime}.$\ Further, let us denote%
\[
H^{\prime}=\left\{  \sum_{\gamma\in F}H_{\gamma}:F\subset\Gamma_{+}\text{ a
finite subset.}\right\}
\]
and%
\[
H^{\prime\prime}=lin\left\{  u_{\gamma}\xi:\xi\in H^{\prime},\gamma\in
\Gamma\right\}  .
\]
Clearly $H^{\prime}$\ is dense in $H_{+}$\ and $H"$\ is dense in $H.$ Define
the linear operator $\widehat{x}$\ on $H"$\ as follows%
\[
\widehat{x}(\sum_{i=1}^{n}u_{\gamma_{i}}\xi_{i})=\sum_{i=1}^{n}u_{\gamma_{i}%
}x\xi_{i}.
\]
We will prove first that $\widehat{x}$\ is well defined. Suppose that
\begin{equation}
\sum_{i=1}^{n}u_{\gamma_{i}}\xi_{i}=0. \label{c4}%
\end{equation}
Since, by Lemma 3.3., $\Gamma$\ is lattice ordered, let $\nu=\inf\left\{
\gamma_{i}:1\leqslant i\leqslant n\right\}  \in\Gamma.$ Since, by hypothesis
$sp(\alpha)=\Gamma,$ we have $\nu\in sp(\alpha).$\ Let $u_{\nu}\in M_{\nu}%
$\ be a unitary operator as in the hypothesis. By multiplying (\ref{c4}) by
$u_{\nu}^{\ast}$\ we get%
\[
\sum_{i=1}^{n}u_{\gamma_{i}-\nu}\xi_{i}=0.
\]
Since $\gamma_{i}-\nu\in\Gamma_{+},i=1,2,...,n$ and $x\in(M_{+})^{\prime}$\ it
follows that%
\[
x\sum_{i=1}^{n}u_{\nu}^{\ast}u_{\gamma_{i}}\xi_{i}=\sum_{i=1}^{n}u_{\nu}%
^{\ast}u_{\gamma_{i}}x\xi_{i}=0.
\]
so%
\[
u_{\nu}\sum_{i=1}^{n}u_{\nu}^{\ast}u_{\gamma_{i}}x\xi_{i}=\sum_{i=1}%
^{n}u_{\gamma_{i}}x\xi_{i}=\widehat{x}(\sum_{i=1}^{n}u_{\gamma_{i}}\xi_{i})=0
\]
so $\widehat{x}$\ is well defined. We will prove next that $\widehat{x}$\ is
continuous. Indeed%
\[
\left\Vert \widehat{x}(\sum_{i=1}^{n}u_{\gamma_{i}}\xi_{i})\right\Vert
=\left\Vert \sum_{i=1}^{n}u_{\gamma_{i}}x\xi_{i}\right\Vert =\left\Vert
u_{\nu}\sum_{i=1}^{n}u_{\nu}^{\ast}u_{\gamma_{i}}x\xi_{i}\right\Vert =
\]%
\[
=\left\Vert \sum_{i=1}^{n}u_{\nu}^{\ast}u_{\gamma_{i}}x\xi_{i}\right\Vert
=\left\Vert x\sum_{i=1}^{n}u_{\nu}^{\ast}u_{\gamma_{i}}\xi_{i}\right\Vert
\leqslant\left\Vert x\right\Vert \left\Vert \sum_{i=1}^{n}u_{\nu}^{\ast
}u_{\gamma_{i}}\xi_{i}\right\Vert =
\]%
\[
\left\Vert x\right\Vert \left\Vert u_{\nu}^{\ast}\sum_{i=1}^{n}u_{\gamma_{i}%
}\xi_{i}\right\Vert =\left\Vert x\right\Vert \left\Vert \sum_{i=1}%
^{n}u_{\gamma_{i}}\xi_{i}\right\Vert .
\]
so $\widehat{x}$\ is continuous. As in the proof of the previous Lemma 3.10.
we can see that $\widehat{x}\in M^{\prime}$\ and $p_{+}\widehat{x}p_{+}=x$, so
$x\in(M^{\prime})_{+}.$

ii) follows from i) by replacing $M$\ with $M^{^{\prime}}.$
\end{proof}

\bigskip

It is worth mentioning that the previous two Lemmas imply, in particular, that
$\left(  M_{+}\right)  "=M_{+}$ but we will not use this fact$.$

\bigskip

In the next two lemmas we will assume that $(M,G,\alpha)\ $is a W*-dynamical
system that satisfies Condition (\textbf{C})\ and $M$ is in standard form.\ We
also assume that $\Gamma$\ is a direct product or sum of archimedean linearly
ordered abelian discrete groups. If $\gamma_{0}\in sp(\alpha)\cap\Gamma
_{+}\smallsetminus\left\{  0\right\}  $and $\gamma\in sp(\alpha)\cap\Gamma
_{+}$ according to Corollary 3.2. there exists $p\in%
\mathbb{Z}
_{+}$ such that either $p\gamma_{0}\leqslant\gamma<(p+1)\gamma_{0}$
or\ $p\gamma_{0}\leqslant\gamma$\ and $\gamma$\ is not comparable with
$(p+1)\gamma_{0}.$ If, in addition, $e_{\gamma}e_{\gamma_{0}}\neq0$ we will
denote for every $\lambda\in%
\mathbb{C}
,\left\vert \lambda\right\vert <1$%
\[
K_{\gamma_{0},\gamma,\lambda}=\overline{\left\{  x(\lambda,\gamma_{0}%
,\gamma,\xi)=\sum_{n\geqslant0}\lambda^{n}u_{\gamma_{0}}^{n}u_{\gamma
-p\gamma_{0}}\xi:\xi\in e_{\gamma_{0}}H_{0}\right\}  }\subset e_{\gamma_{0}%
}H_{+}.
\]
and%
\[
L_{\gamma_{0},\gamma,\lambda}=\overline{\left\{  y(\lambda,\gamma_{0}%
,\gamma,\xi)=\sum_{n\geqslant0}\lambda^{n}w_{\gamma_{0}}^{n}w_{\gamma
-p\gamma_{0}}\xi:\xi\in e_{\gamma_{0}}H_{0}\right\}  }\subset e_{\gamma_{0}}H
\]

\bigskip

\textbf{Lemma 3.12. }\textit{Let}\textbf{\ }$(M,G,\alpha)\ $\textit{be a
W*-dynamical system that satisfies Condition (\textbf{C}). }$\Gamma$\textit{
and }$\gamma_{0},\gamma\in sp(\alpha)\cap\Gamma_{+}$ \textit{be as
above}.\ \textit{Then, }%
\[
K_{\gamma_{0}}=lin\left\{  K_{\gamma_{0},\gamma,\lambda}:\gamma\in
sp(\alpha)\cap\Gamma_{+},\lambda\in%
\mathbb{C}
,\left\vert \lambda\right\vert <1\right\}
\]
\textit{ is dense in }$e_{\gamma_{0}}H_{+}.$

\bigskip

\begin{proof}
Notice first that if $\gamma_{0},\gamma$\ are as in the hypothesis of the
lemma, then, by the definition of $e_{\gamma_{0}}$ in the Condition
(\textbf{C})$,$we have that $p\gamma_{0}\in sp(\alpha)$\ for every $p\in%
\mathbb{Z}
_{+}$\ and, by Lemma 3.6., $\gamma-p\gamma_{0}\in sp(\alpha),$\ so
$u_{\gamma-p\gamma_{0}}$ exists, and thus the definition of $K_{\gamma
_{0},\gamma,\lambda}$\ in the hypothesis of the Lemma is consistent. Now,
taking $\lambda=0$\ in $K_{\gamma_{0},\gamma,\lambda},$\ it follows that
$e_{\gamma_{0}}u_{\gamma-p\gamma_{0}}H_{0}\subset K_{\gamma_{0},\gamma}\subset
K_{\gamma_{0}}.$ In particular, for $p=0\ $(so, when either $0\leqslant
\gamma<\gamma_{0}$ or $0\leqslant\gamma\ $and $\gamma\ $is not comparable with
$\gamma_{0}$),\ we have $e_{\gamma_{0}}u_{\gamma}H_{0}=e_{\gamma_{0}}%
H_{\gamma}\subset K_{\gamma_{0}}$. We will prove that $e_{\gamma_{0}}%
H_{\gamma}=u_{\gamma}e_{\gamma_{0}}H_{0}\subset K_{\gamma_{0}}$ for every
$\gamma\in sp(\alpha)\cap\Gamma_{+}.$ This fact will imply that $\sum
_{\gamma\in\Gamma_{+}}e_{\gamma_{0}}H_{\gamma}\subset K_{\gamma_{0}}$ and
therefore, $e_{\gamma_{0}}H_{+}=\overline{\sum_{\gamma\in\Gamma_{+}}%
e_{\gamma_{0}}H_{\gamma}}\subset\overline{K}_{\gamma_{0}},$ so $K_{\gamma_{0}%
}$ is dense in $e_{\gamma_{0}}H_{+}$ as claimed. We will prove first that
$e_{\gamma_{0}}u_{\gamma_{0}}^{p}u_{\gamma-p\gamma_{0}}H_{0}\subset
K_{\gamma_{0},\gamma}.$ The case $p=0$\ was proved above. Suppose that
$p>0.$\ We will prove by induction on $k$\ that%
\[
e_{\gamma_{0}}u_{\gamma_{0}}^{k}u_{\gamma-p\gamma_{0}}H_{0}\subset
K_{\gamma_{0}}%
\]
for every $k,$\ in particular for $k=p.$\ If $k=0,$\ the above inclusion
follows, as noticed at the beginning of this proof from the definition of
$K_{\gamma_{0},\gamma}$\ for $\lambda=0.$ Suppose by induction that
$u_{\gamma_{0}}^{l}u_{\gamma-p\gamma_{0}}\xi\in K_{\gamma_{0}}$ for
$l=0,1,...k-1.$Then%
\[
\sum_{n\geqslant k}\lambda^{n}u_{\gamma_{0}}^{n}u_{\gamma-p\gamma_{0}}\xi\in
K_{\gamma_{0}}%
\]
for every $\xi\in e_{\gamma_{0}}H_{0},\lambda\in%
\mathbb{C}
,\left\vert \lambda\right\vert <1.$\ Thus%
\[
\lambda^{k}u_{\gamma_{0}}^{k}u_{\gamma-p\gamma_{0}}\xi+\sum_{n\geqslant
k+1}\lambda^{n}u_{\gamma_{0}}^{n}u_{\gamma-p\gamma_{0}}\xi\in K_{\gamma_{0}}%
\]
By dividing the above relation by $\lambda^{k}$, $\lambda\neq0$ and then
taking the limit as $\lambda\rightarrow0,$ we get that $u_{\gamma_{0}}%
^{k}u_{\gamma-p\gamma_{0}}\xi\in K_{\gamma_{0}},$ so, in particular,
$u_{\gamma_{0}}^{p}u_{\gamma-p\gamma_{0}}\xi\in K_{\gamma_{0}}$ for every
$\xi\in H_{0}$. Since, obviously, $u_{\gamma_{0}}^{p}(u_{\gamma_{0}}^{\ast
})^{p}=e_{\gamma_{0}}$ for $p>0$ and, by Lemma 2.5. ii), $M_{-p\gamma_{0}%
}M_{\gamma}\subset M_{\gamma-p\gamma_{0}},$ it follows that
\[
e_{\gamma_{0}}H_{\gamma}=e_{\gamma_{0}}u_{\gamma}H_{0}=u_{\gamma_{0}}%
^{p}(u_{\gamma_{0}}^{\ast})^{p}u_{\gamma}H_{0}\subset u_{\gamma_{0}}%
^{p}u_{\gamma-p\gamma_{0}}H_{0}\subset K_{\gamma_{0}}%
\]
and we are done.
\end{proof}

\bigskip

The following lemma can be proven similarly with the previous Lemma 3.12.

\bigskip

\textbf{Lemma 3.13. }\textit{Let}\textbf{\ }$(M,G,\alpha)\ $\textit{be a
W*-dynamical system that satisfies Condition (\textbf{C}), }$\Gamma$\textit{
and }$\gamma_{0},\gamma\in sp(\alpha)\cap\Gamma_{+}$ \textit{be as
above}.\ \textit{Then, }%
\[
L_{\gamma_{0}}=lin\left\{  L_{\gamma_{0},\gamma,\lambda}:p\gamma_{0}%
\leqslant\gamma<(p+1)\gamma_{0}\text{ \ \textit{for some} }p\in%
\mathbb{Z}
_{+},\lambda\in%
\mathbb{C}
,\left\vert \lambda\right\vert <1\right\}
\]
\textit{ is dense in }$e_{\gamma_{0}}H_{+}.$

\bigskip

\begin{proof}
\textbf{of Theorem 3.8. }We\textbf{ }will prove first that $M_{+}\subset
B(H_{+})$ is reflexive and then apply Lemma 3.4. and the subsequent discussion
to infer that $M_{+}$\ is hereditarily reflexive. Let $\gamma_{0}\in
sp(\alpha)\cap\Gamma_{+}.$ Since $u_{\gamma_{0}}\in M,$\textbf{ }it follows
that\textbf{ }$u_{\gamma_{0}}f=fu_{\gamma_{0}}$\ for every projection
$f\in(M^{\prime}),$ in particular for every projection $f\in(M^{\prime})_{0}%
.$\ By\textbf{ }Lemma 2.8. iv), since $\overline{(M^{^{\prime}})_{\gamma}%
H_{0}}=H_{\gamma},\gamma\in sp(\alpha),$\ we have, in particular that
$(M^{^{\prime}})_{0}p_{+}\subset B(H_{+}).$ Therefore, for every projection
$f\in(M^{^{\prime}})_{0},\ fH_{+}$belongs to $Lat(M_{+}).$ Let $x\in
algLat(M_{+})\subset B(H_{+}).$ We will prove that $x\in((M^{^{\prime}}%
)_{+})^{^{\prime}}$\ and, then, applying Lemma 3.11. ii) it will follow that
$x\in M_{+},$\ so $M_{+}$\ is a reflexive operator algebra. The way to prove
this fact is to use Lemma 3.14. to show that $x^{\ast}\in((M^{^{\prime}}%
)_{+}^{\ast})^{^{\prime}}$\ and then, clearly, it will follow that
$x\in((M^{^{\prime}})_{+})^{^{\prime}}=M_{+}.$As noticed above, $fH_{+}%
\in\mathit{LatM}_{+}$ for every projection $f\in(M^{^{\prime}})_{0},$so
$xf=fx$\ and therefore $x^{\ast}f=fx^{\ast}$ for every projection
$f\in(M^{^{\prime}})_{0}.$ It follows that $x$\ commutes with every element of
$(M^{^{\prime}})_{0}.$ We will prove next that $xw_{\gamma}=w_{\gamma}x$ for
every $\gamma\in sp(\alpha^{\prime})\cap\Gamma_{+}=sp(\alpha)\cap\Gamma_{+}$
and then apply Lemma 3.10. ii) to infer that $x\in((M^{^{\prime}}%
)_{+})^{^{\prime}}=M_{+}.$ To this end, let $\gamma_{0}\in sp(\alpha
)\cap\Gamma_{+}.$\ If $\gamma_{0}=0,$ then as convened, $u_{\gamma_{0}%
}=w_{\gamma_{0}}=I,$\ so nothing to prove. Let $\gamma_{0}\in sp(\alpha
)=sp(\alpha^{\prime}),\gamma_{0}>0.$ Denote by $T_{u_{\gamma_{0}}}%
$\ (respectively $T_{w_{\gamma_{0}}}$) the operator $u_{\gamma_{0}}\in M_{+}$
(respectively $w_{\gamma_{0}}\in(M^{\prime})_{+}$) defined on $H_{+}.$ Then
the adjoints of $T_{u_{\gamma_{0}}},T_{w_{\gamma_{0}}}$ on $H_{+}$ are%
\[
T_{u_{\gamma_{0}}}^{\ast}\xi_{\gamma}=0\text{ if }0\leqslant\gamma<\gamma
_{0}\text{ and }T_{u_{\gamma_{0}}}^{\ast}\xi_{\gamma}=u_{\gamma_{0}}^{\ast}%
\xi_{\gamma}\text{ if }\gamma\geqslant\gamma_{0}.
\]
and similarly%
\[
T_{w_{\gamma_{0}}}^{\ast}\xi_{\gamma}=0\text{ if }0\leqslant\gamma<\gamma
_{0}\text{ and }T_{w_{\gamma_{0}}}^{\ast}\xi_{\gamma}=w_{\gamma_{0}}^{\ast}%
\xi_{\gamma}\text{ if }\gamma\geqslant\gamma_{0}.
\]
Since $u_{\gamma_{0}},w_{\gamma_{0}}$\ (so $T_{u_{\gamma_{0}}},T_{w_{\gamma
_{0}}})$ commute, it follows that $T_{u_{\gamma_{0}}}^{\ast},T_{w_{\gamma_{0}%
}}^{\ast}$ commute as well. Notice that if, for $\lambda\in%
\mathbb{C}
,\left\vert \lambda\right\vert <1,$ we denote
\[
\widetilde{L}_{\gamma_{0,\lambda}}=\left\{  \xi\in H_{+}:T_{w_{\gamma_{0}}%
}^{\ast}\xi=\lambda\xi\right\}  .
\]
then, since $T_{w_{\gamma_{0}}}^{\ast}$ commutes with $(M_{+})^{\ast},$\ we
have $\widetilde{L}_{\gamma_{0},\gamma,\lambda}\in Lat(M_{+})^{\ast}$ and,
since $x\in algLat(M_{+})$ it follows that$.L_{\gamma_{0},\gamma,\lambda}\in
Lat(x^{\ast})$ for every $\lambda\in%
\mathbb{C}
,\left\vert \lambda\right\vert <1.$ Since $\widetilde{L}_{\gamma_{0}%
,\gamma,\lambda}$\ consists of eigenvectors of $T_{w_{\gamma_{0}}}^{\ast}%
$\ for the eigenvalue $\lambda,$\ then, it follows that $x^{\ast}%
T_{w_{\gamma_{0}}}^{\ast}\xi=T_{w_{\gamma_{0}}}^{\ast}x^{\ast}\xi$ for every
$\xi\in\widetilde{L}_{\gamma_{0,\lambda}}.$\ On the other hand, it is clear
that
\[
L_{\gamma_{0},\gamma,\lambda}\subset\widetilde{L}_{\gamma_{0,\lambda}}%
\]
where $L_{\gamma_{0},\gamma,\lambda}$ is as in Lemma 3.13., so $x^{\ast
}T_{w_{\gamma_{0}}}^{\ast}\xi=T_{w_{\gamma_{0}}}^{\ast}x^{\ast}\xi$ for every
$\xi\in L_{\gamma_{0,\gamma,\lambda}}.$ By Lemma 3.13.,
\[
L_{\gamma_{0}}=lin\left\{  L_{\gamma_{0},\gamma,\lambda}:p\gamma_{0}%
\leqslant\gamma<(p+1)\gamma_{0}\text{ \ \textit{for some} }p\in%
\mathbb{Z}
_{+},\lambda\in%
\mathbb{C}
,\left\vert \lambda\right\vert <1\right\}
\]
\ is dense in $e_{\gamma_{0}}H_{+},$ and therefore $x^{\ast}T_{w_{\gamma_{0}}%
}^{\ast}\xi=T_{w_{\gamma_{0}}}^{\ast}x^{\ast}\xi$ for every $\xi\in
e_{\gamma_{0}}H_{+}$, so $x^{\ast}$ commutes with $T_{w_{\gamma_{0}}}^{\ast
},\gamma_{0}\in sp(\alpha)\cap\Gamma_{+}\ $and therefore with $\left(
M^{\prime}\right)  _{+}^{\ast}.$ Since $x^{\ast}$ commutes with $\left(
M^{\prime}\right)  _{+}^{\ast},$ it follows that $x$ commutes with $\left(
M^{\prime}\right)  _{+}$, so by Lemma 3.10. ii) $x\in M_{+}$\ so $M_{+}$\ is
reflexive. Finally, by applying Lemma 3.4. and the discussion following it, we
see that $M_{+}$\ is hereditarily reflexive and we are done.
\end{proof}

\bigskip

Corollary 3.14. below\ refers to the Example a) to Condition (\textbf{C}).

\bigskip

\textbf{Corollary 3.14. }\textit{Let }$(M,G,\alpha)$ \textit{be a
W}*\textit{-dynamical system with }$G$\textit{\ compact abelian and }%
$M$\textit{\ a finite W*-algebra in standard form such that }$Z(M_{0})\subset
Z(M)$\textit{. Suppose that the dual }$\Gamma$\textit{\ of }$G$\textit{\ has
an archimedean linear order. Then }$M_{+}\subset B(H_{+})$ \textit{is
reflexive.}

\bigskip

\begin{proof}
According to [18, Theorem 2.3.], if $M$\ is finite and $Z(M_{0})\subset Z(M)$,
then the Condition (\textbf{C}) is satisfied and therefore the result follows
from Theorem 3.8.
\end{proof}

\bigskip

\begin{proof}
of Theorem 3.9. The proof is very similar with the proof of Theorem 3.8. The
only modification is using Lemma 3.11 \ instead of Lemma 3.10.
\end{proof}

\bigskip

The next Corollary refers to Examples b) to Condition (\textbf{C}).

\bigskip

\textbf{Corollary 3.15. }\textit{Let} $(M,G,\alpha)$ \textit{be a
W}*\textit{-dynamical system with }$M$\textit{\ an injective von Neumann
algebra in standard form and }$G$\ \textit{a compact abelian group}
\textit{such that the dual }$\Gamma$\textit{\ of }$G$\textit{\ is a direct
product (or a direct sum) of archimedean linearly ordered discrete
groups.\ Suppose that }$\alpha$\textbf{ }\textit{is prime and faithful and it
is not the case that }$M$\textit{\ is of type }$\mathit{III}$\textit{ and
}$M_{0}$\textit{\ is of type }$\mathit{II}_{1}$\textit{. Then }$M_{+}$\textit{
is hereditarily reflexive.}

\bigskip

\begin{proof}
Since $\alpha$\ is faithful, we have $sp(\alpha)=\Gamma$\ (see for instance
[21, Lemma 2.2.]). By [21, Theorem 2.3.] each spectral subspace $M_{\gamma}%
$\ contains a unitary operator.\ The conclusion of the Corollary follows from
Theorem 3.9.
\end{proof}

\bigskip

The concept of nonselfadjoint crossed product, or more generally that of
w*-semicrossed product were defined in [3], [11], [16].

\bigskip

\textbf{Corollary 3.16. }\textit{Let }$(M_{0},\Gamma,\beta)$\textit{\ be a
W*-dynamical system such that }$M_{0}\subset B(H_{0})$\textit{\ is in standard
form and }$\Gamma$\textit{\ is a discrete abelian group. Suppose that }%
$\Gamma\subseteq\Pi\Gamma_{\iota}$ \textit{is a direct product (or a direct
sum) of archimedean linearly ordered groups. Let} $M=M_{0}\times_{\alpha
}\Gamma\subset B(l^{2}(\Gamma,H_{0}))$ be \textit{the corresponding crossed
product}$\mathit{.}$\textit{ Then, the non selfadjoint crossed product }%
$M_{+}=M_{0}\times_{\beta}\Gamma_{+}$ \textit{(i.e. the algebra of elements of
}$M$\textit{\ with non negative spectrum) is a hereditarily reflexive operator
algebra in }$B(H_{+})$\textit{\ where }$H_{+}=l^{2}(\Gamma_{+},H_{0}).$

\bigskip

\begin{proof}
If $G$\ denotes the (compact) dual of $\Gamma,$\ and $\alpha=\widehat{\beta}%
$\ is the dual action of $\beta$\ on $M,$\ then, consider the canonical
conditional expectation $P_{0}:M\rightarrow M_{0}.$\ Since $M_{0}\subset
B(H_{0})$\ is in standard form, from Lemma 2.1. and the subsequent discussion
it follows that $M\subset B(l^{2}(\Gamma,H_{0}))\ $is in standard form and by
the definition of the crossed product, $M_{\gamma}$\ contains a unitary
operator $u_{\gamma}$\ for every $\gamma\in\Gamma$\ and thus the W*-dynamical
system $(M,G,\widehat{\beta})$ satisfies the conditions of Corollary 3.15.
\end{proof}

\bigskip

\textbf{Corollary 3.17. }\textit{Let }$(M_{0},\Gamma,\beta)$\textit{\ be a
W*-dynamical system such that }$M_{0}\subset B(H_{0})$\textit{\ is a maximal
abelian von Neumann algebra and }$\Gamma$\textit{\ is a discrete abelian
group. Suppose that }$\Gamma\subseteq\Pi\Gamma_{\iota}$ \textit{is a direct
product (or a direct sum) of archimedean linearly ordered groups. Let}
$M=M_{0}\times_{\alpha}\Gamma\subset B(l^{2}(\Gamma,H_{0}))$ be \textit{the
corresponding crossed product}$\mathit{.}$\textit{ Then, the non selfadjoint
crossed product }$M_{+}=M_{0}\times_{\beta}\Gamma_{+}$ \textit{ is a
hereditarily reflexive operator algebra in }$B(H_{+})$\textit{\ where }%
$H_{+}=l^{2}(\Gamma_{+},H_{0}).$

\bigskip

\begin{proof}
Since $M_{0}\subset B(H_{0})$\textit{\ }is a maximal abelian von Neumann
algebra, by Proposition 2.4. it is spatially isomorphic with its standard
form, so the result folows from the previous Corollary 3.16.
\end{proof}

\bigskip

In [2, Corollary 5.14.] it is stated that if $M_{0}\subset B(H_{0})$\ is a
maximal abelian von Neumann algebra and $\Gamma=%
\mathbb{Z}
^{d},d\in%
\mathbb{N}
,$\ then $M_{0}\times_{\beta}\Gamma_{+}$\ is reflexive, so the Corollary 3.18.
above extends that result by showing also hereditary reflexivity in the
special case $\Gamma=%
\mathbb{Z}
^{d},d\in%
\mathbb{N}
.$

\bigskip

\ \ \ \ \ \ \ \ \ \ \ \ \ \ \ \ \ \ \ \ \ \ \ \ \ \ \ \ \textbf{REFERENCES}

1. W. B. Arveson, The harmonic analysis of automorphism groups, operator
algebras and applications, Proc. Sympos. Pure Math., vol. 38, Amer. Math.
Soc., Providence, R. I., 1982.

2. R. T. Bickerton, and E. T. A. Kakariadis, Free multivariate w*-semicrossed
products: reflexivity and the bicommutant property, Canad. J. Math. 70(2018), 1201--1235.

3. K, R. Davidson, A. H. Fuller and E. T. A. Kakariadis, Semicrossed Products
of Operator Algebras by Semigroups, Memoirs of the AMS, Vol. 247, 2017.

4. D. W. Hadwin and E.A. Nordgren, Subalgebras of reflexive algebras. J.
Operator Theory 7 (1982), 3--23.

5. L. Helmer, Reflexivity of non-commutative Hardy algebras, J. Funct. Anal.
272(2017), 2752--2794.

6. V. F. R. Jones, Prime actions of compact abelian groups on the hyperfinite
type II$_{1}$ factor, J. Operator Theory, 9(1983), 181-186.

7. V. F. R. Jones and M. Takesaki, Actions of compact abelian groups on
semifinite injective factors, Acta Math. 153 (1984), 213--258.

8. R. V. Kadison and J. R. Ringrose, Fundamentals of the theory of operator
algebras, Vol. II Advanced Theory,Academic Press 1986.

9. E. T. A. Kakariadis, Semicrossed products and reflexivity, J. Operator
Theory, 67(2012), 379-395.

10. A I Loginov and V S \v{S}ul'man, Hereditary and intermediate reflexivity
of W*-algebras, Mathematics of the USSR-Izvestiya 9(1975), 1189-1202.

11. M. McAsey, P. S. Muhly, and K.-S. Saito, Nonselfadjoint crossed products
(invariant subspaces and maximality), Trans. Amer. Math. Soc., 248(1979), 381--409.

12. P. S. Muhly and B. Solel, Hardy algebras, W*-correspondences and
interpolation theory, Math. Ann. 330 (2004), 353--415.

13. G. K. Pedersen, C*-algebras and their automorphism groups, Academic Press 1979.

14. C. Peligrad, Reflexive operator algebras on noncommutative Hardy spaces,
Math. Annalen, 253(1980), 165-175.

15. C. Peligrad, Invariant subspaces of algebras of analytic elements
associated with periodic flows on von neumann algebras, Houston J. Math.,
42(2016), 1331-1445.

16. J. R. Peters, Semicrossed products of C*-algebras. J. Funct. Anal.,
59(1984), 498--534.

17. H. Radjavi and P. Rosenthal,\ Invariant subspaces, 2nd edition, Dover
Publications, Mineola, New, York, 2003\textit{.}

18. K.-S. Saito, Nonselfadjoint subalgebras associated with compact abelian
group actions on finite von Neumann algebras, Tohoku Math. Journ. 34(1982), 485-494.

19. D. Sarason, Invariant subspaces and unstarred operator algebras, Pacific
J. Math., 17(1966), 511-517.

20. S. Stratila, Modular theory in operator algebras, Bucharest; Abacus Press,
Tunbridge Wells, 1981

21. K. Thomsen, Compact abelian prime actions on von Neumann algebras, Trans.
Amer. Math. Soc., 315(1989), 255-273.

\bigskip
\end{document}